\documentclass[10pt]{amsart}
\usepackage{amsmath,amsfonts,amsthm,bbm, amssymb}
\usepackage[cmtip, all]{xy}
\usepackage[left=3cm, right=3cm, top = 2.5cm, bottom = 2.5cm ]{geometry}
\usepackage{comment}
 
\usepackage{enumerate}
\usepackage[usenames]{xcolor}
\usepackage{comment}
\usepackage{graphicx}
\usepackage{amssymb}
\usepackage{stmaryrd}
\usepackage{dsfont, mathrsfs,bold-extra}
\usepackage{setspace}
\usepackage[color, notref, notcite, final]{showkeys}     
\definecolor{refkey}{gray}{.5}   % graylevel for refs
\definecolor{labelkey}{gray}{.5} % graylevel for labels
\definecolor{Red}{rgb}{1,0,0}

\usepackage{soul}

\newtheorem{thmintro}{Theorem}
\newtheorem{corintro}[thmintro]{Corollary}

\newtheorem{thm}{Theorem}[section]
\newtheorem{prop}[thm]{Proposition}
\newtheorem{lem}[thm]{Lemma} 
\newtheorem{cor}[thm]{Corollary}

\newtheorem{claim}[thm]{Claim}

\theoremstyle{definition}

\newtheorem{df}[thm]{Definition}
\theoremstyle{remark}
\newtheorem{rmk}[thm]{Remark}
\newtheorem{remks}[thm]{Remarks}

\newtheorem{exm}[thm]{Example}
\newtheorem{exms}[thm]{Examples}
\newtheorem{notat}[thm]{Notation}
\numberwithin{equation}{section}

{\hfill$\square$\end{defn}}
{\hfill$\square$\end{remk}}
{\hfill$\square$\end{remks}}
{\hfill$\square$\end{exm}}
{\hfill$\square$\end{exms}}
    \newcommand{\bN}{\mathbb{N}}
    
	\newcommand{\Z}{\mathbb{Z}}

	\newcommand{\Q}{\mathbb{Q}}
	
	\renewcommand{\P}{\mathbb{P}}

    \def\Zbk{\ol{Z}_k}
    \providecommand{\cub}[1]{\ol{\square}^{#1}}
\providecommand{\frac}[1]{\operatorname{Frac}(#1)}  %Fraction field
\providecommand{\Spec}[1]{\operatorname{Spec}(#1)} 	%
\providecommand{\Proj}[1]{\operatorname{Proj}(#1)}  %Proj functor - 2 versions (the capital one is the usual one. The other one is for backw. compatibility)
\newcommand{\id}{\operatorname{id}}   			%identity
\renewcommand{\H}{\operatorname{H}}				     % \H for homological operators
\newcommand{\CH}{\operatorname{CH}}					 %  Chow group/Chow ring
\renewcommand{\lim}{\operatorname{lim}}			 %limit  - I don't like the standard one

	\def\Cb{\ol{C}}
	\newcommand{\tensor}{\otimes}
	\renewcommand{\cong}{\simeq}
	
	\newcommand{\ol}{\overline}

\def\Xb{\ol{X}} %
\def\Yb{\ol{Y}}
\def\Wb{\ol{W}}
\def\Vb{\ol{V}}

\newcommand{\cC}{\mathcal{C}}

\newcommand{\cF}{\mathcal{F}}

\newcommand{\cO}{\mathcal{O}}
\newcommand{\cP}{\mathcal{P}}

\newcommand{\cX}{\mathcal{X}}

\renewcommand{\epsilon}{\varepsilon}
\renewcommand{\phi}{\varphi}

\newcommand{\cubb}{\square}

\def\RSC{\mathbf{RSC}}

\makeatletter
\newcommand{\colim@}[2]{%
  \vtop{\m@th\ialign{##\cr
    \hfil$#1\operator@font hocolim$\hfil\cr
    \noalign{\nointerlineskip\kern1.5\ex@}#2\cr
    \noalign{\nointerlineskip\kern-\ex@}\cr}}%
}
\newcommand{\hocolim}{%
  \mathop{\mathpalette\colim@{\rightarrowfill@\textstyle}}\nmlimits@
}
\makeatother

\makeatletter
\newcommand{\ncolim@}[2]{%
  \vtop{\m@th\ialign{##\cr
    \hfil$#1\operator@font colim$\hfil\cr
    \noalign{\nointerlineskip\kern1.5\ex@}#2\cr
    \noalign{\nointerlineskip\kern-\ex@}\cr}}%
}
\newcommand{\varcolim}{%
  \mathop{\mathpalette\ncolim@{\rightarrowfill@\textstyle}}\nmlimits@
}
\makeatother

\makeatletter
\newcommand{\varvarcolim}{%
  \mathop{\mathpalette\ncolim@{}}\nmlimits@
}
\makeatother

\makeatletter
\newcommand{\holim@}[2]{%
  \vtop{\m@th\ialign{##\cr
    \hfil$#1\operator@font holim$\hfil\cr
    \noalign{\nointerlineskip\kern1.5\ex@}#2\cr
    \noalign{\nointerlineskip\kern-\ex@}\cr}}%
}
\newcommand{\holim}{%
  \mathop{\mathpalette\holim@{\leftarrowfill@\textstyle}}\nmlimits@
}
\makeatother

\def\hocof{\operatorname{hocof}}
\def\hofib{\operatorname{hofib}}
\definecolor{winered}{rgb}{0.8,0,0}
\usepackage[pdfauthor={Federico Binda}, %
				pdftitle={MilnorChowMod} ,colorlinks=true]{hyperref}
\hypersetup{
	bookmarksnumbered=true,
	linkcolor={winered},
	citecolor=blue,
}

\usepackage{float}
\newcounter{elno}   
\newenvironment{romanlist}{
                         \begin{list}{\roman{elno})
                                     }{\usecounter{elno}}
                      }{
                         \end{list}}
                         
\newcounter{elno-abc}   

\begin{document}
\newcommand{\nocontentsline}[3]{}
\newcommand{\tocless}[2]{\bgroup\let\addcontentsline=\nocontentsline#1{#2}\egroup}

\author{Federico Binda}
\title{A cycle class map from Chow groups with modulus to relative $K$-theory}

\AtEndDocument{\bigskip{\footnotesize%
  (F.~Binda) \textsc{Fakult\"at f\"ur Mathematik,  Universit\"at Regensburg, 93040 Regensburg, Germany.} 
  \textit{E-mail address},  \texttt{federico.binda@ur.de}}}
\thanks{Supported by the DFG SBF/CRC 1085 ``Higher Invariants''.}
\maketitle
\date{\today}
%%%%%abstract
\renewcommand{\div}{\operatorname{div}}
\begin{abstract}
Let $\ol{X}$ be a smooth quasi-projective $d$-dimensional variety over a field $k$ and let $D$ be an effective, non-reduced, Cartier divisor on it such that its support is strict normal crossing. In this note, we construct cycle class maps from (a variant of) the higher Chow group with modulus of the pair $(\ol{X};D)$ in the range $(d+n, n)$ to the relative $K$-groups $K_n(\ol{X}; D)$  for every $n\geq 0$. 
\end{abstract}
\tableofcontents

\section*{Introduction}
\tocless\subsection{}In \cite{Bloch86}, Bloch introduced his higher Chow groups $\CH^q( X,p)$ for a variety $X$ over a field $k$ with the goal of defining an integral cohomology theory which rationally gives the weight-graded pieces of Quillen's algebraic $K$-theory $K_p(X)^{(q)}$. The higher Chow groups are defined as the homology of a certain explicit complex of algebraic cycles, and when $X$ is a smooth quasi-projective variety over $k$, they are known to agree with the motivic cohomology groups defined by Voevodsky \cite{VSF} (see e.g. \cite[Theorem 1.2]{LevineLoc}). To compare the higher Chow groups of $X$ and the relevant graded pieces of $K_p(X)$,  Bloch constructed functorial Chern classes 
\[c_{q,p}\colon K_{2q-p}(X)\to \CH^q(X, 2q-p)\]
which can be shown to satisfy standard properties.

The situation is much more mysterious when $X$ is not smooth. The theory of algebraic cycles with modulus,  developed in \cite{BS} and motivated by the work of Kerz and Saito \cite{KS1}, is an evolution of the theory of additive Chow groups introduced by Bloch-Esnault \cite{BEAdditive} (and later developed by Park, Krishna-Levine, Krishna-Park) with the aim of understanding the algebraic $K$-theory of non-reduced schemes, such as the truncated polynomial ring $k[t]/(t^n)$, in terms of algebraic cycles. In this sense, it can be seen as an attempt to generalize Bloch's higher Chow groups in order to capture some extra infinitesimal information.

\tocless\subsection{}To fix the ideas, let $\ol{X}$ be a scheme, equidimensional and of finite type over $k$, and let $D$ be an effective Cartier divisor on it. The higher Chow groups with modulus are defined as the homology of the relative cycle complex $z^q(\ol{X}|D, *)$, a subcomplex of Bloch's cycle complex  $z^q(\ol{X}, *)$ consisting of cycles satisfying an additional condition along $D$ (see \ref{def:definitionCrXD-both-MC} for the precise definition). It is expected that (a sheafified version of) the cycle complex with modulus $z^q(\ol{X}|D, *)[-2q]$ of the pair $(\ol{X};D)$ can play the role of Bloch's complex in the framework of a relative motivic theory (or theory of \textit{motives with modulus}), as currently under development  in the program of Kahn-Saito-Yamazaki (see \cite{KSY2}, \cite{KSY}, \cite{KSY-RecAndMotives}). From this perspective, there should be a close relationship between the groups of relative cycles and the relative $K$-groups $K_*(\ol{X};D)$, defined as the homotopy groups of the homotopy fiber of the restriction map $\mathbf{K}(\ol{X})\to \mathbf{K}(D)$. The problem of comparing relative $K$-groups with Chow groups with modulus (or additive Chow groups) is central in  works of many authors (we recall here  \cite{BK}, \cite{BEAdditive},  \cite{KrishnaOnCycles}, \cite{KL}, \cite{ParkAMJ}, \cite{Rul},  \cite{RulSaito}). 

An approach conceptually similar to the one originally proposed by Bloch (although significant new ingredients are required) has been recently exploited by Iwasa-Kai, culminating with the construction of \textit{relative Chern classes} (see \cite[Theorem 5.8]{IwasaKai})
\[c_{q,p}\colon K_{2q-p}(\ol{X};D) \to \CH_{\rm Nis}^q(\ol{X}|D, 2q-p), \]  
where $\CH_{\rm Nis}^q(\ol{X}|D, 2q-p)$ denotes the Nisnevich hypercohomology group $\mathbb{H}^{*}(\ol{X}_{\rm Nis}, z^q(-|-\times D, *))$ (i.e. the relative motivic cohomology of the pair as introduced in \cite{BS}). The classes $c_{q,p}$ are functorial in $(\ol{X};D)$ and coincide with Bloch's Chern classes when $D=\emptyset$.

\tocless \subsection{} In the other direction, not much progress has been made towards the construction of cycle class maps \textit{from} the groups of cycles with modulus \textit{to} the relative $K$-groups. One of the reasons for this is the intrinsic difficulty in showing that \textit{a priori} defined relations (such as the one given by the modulus condition) among cycles give rise to relations in the $K$-groups, and the fact that one is naturally led to consider $K$-groups of singular, non-reduced schemes in doing so. We investigated the already challenging relationship between the group $\CH_0(\ol{X}|D) = \CH_0(\ol{X}|D, 0)$ of zero-cycles with modulus and the relative $K_0$ group $K_0(\ol{X};D)$ in \cite{BK}.

A standard method, which dates back to the works of Bloch and Levine, to construct classes in higher $K$-groups $K_n(X)$ of a regular $k$-variety $X$ is to construct classes in a suitable multi-relative $K_0$ group, namely
\[K_0(X\times \square^n_{\bullet}; X\times \partial \square^n_{\bullet})\]
where $\square^n \cong \mathbb{A}^n_k$ and $\partial \square^n$ is the standard boundary divisor on $\square^n$ (see \ref{sec:rel-cycle-complexes}). The homotopy invariance property of $K$-groups of regular schemes gives in fact a natural isomorphism 
\[K_0(X\times \square^n; X\times \partial \square^n) \xrightarrow{\simeq} K_n(X),\]
 compatible with Adams operations on both sides. It is then possible to show that a codimension $p$ cycle in $X\times \square^n$, in good position with respect to all the faces of $\square^n$, gives rise to a class in $K_0(X\times \square^n; X\times \partial \square^n)^{(p)}$, and therefore in $K_n(X)^{(p)}$ (see \cite[Lemma 2.2 and Theorem 3.1]{LevineBlochrevisited}). 

Unfortunately (actually, this is one of the main features of the theory), the relative $K$-groups $K_n(\ol{X};D)$ are typically far from being $\mathbb{A}^1$-invariant, so that the classical argument cannot go through. Our approach in this work is to replace the homotopy invariance property with an appropriate application of the $\mathbb{P}^1$-bundle formula, available without regularity assumptions. This allows us to construct a new model for the relative $K$-theory of a pair $(\ol{X};D)$ (see Proposition \ref{lem:delooping}).

To some extent, this can be seen as $K$-theoretic analogue of the $\ol{\square}^{(-1)}=(\mathbb{P}^1, -\infty)$-invariance property of higher Chow groups with modulus, recently established by Miyazaki \cite{Miyazaki}, and it's the starting point for the construction of a motivic homotopy category for modulus \textit{triples}, without $\mathbb{A}^1$-invariance, which is the content of \cite{BindaMotModulus} (see also \cite[Chapter II]{BThesis}).

\tocless\subsection{} We can now state our main Theorem.  
\begin{thmintro}[see Theorem \ref{thm:main-theorem-cycle-class}]\label{thm:thm-intro} Let $\ol{X}$ be a smooth quasi-projective $k$-variety of dimension $d$ and let $D$ be an effective Cartier divisor on it. Assume that the support of $D_{\rm red}$ is a strict normal crossing divisor. Then the cycle class map \eqref{eq:def-of-cycleclass-map}, obtained by sending a zero cycle $\alpha \in z^{d+n}(\ol{X}|D, n)$ to its fundamental class (see \ref{sssec:fund-class}),  induces a group homomorphism
    \[{\rm cyc}^{d+n}_{\ol{X}|D}\colon \CH^{d+n}(\ol{X}|D, n)_{M_{\text{ssup}}, \Q} \to K_n(\ol{X};D)^{(d+n)}\]
    for each $n\geq 0$, where $K_n(\ol{X};D)^{(d+n)}$ is the $k^{d+n}$ characteristic subspace of the Adams operation $\psi^k$  acting on $K_n(\ol{X};D)_{\Q}$ (for a fixed integer $k>1$, see \ref{ssec:Adams}  for details). When $D=\emptyset$, the maps ${\rm cyc}^{d+n}_{\ol{X}|\emptyset}$ agree with the cycle class maps defined by Bloch \cite{Bloch86} and Levine \cite{LevineBlochrevisited}. 
    \end{thmintro}
   % The maps  ${\rm cyc}^{d+n}_{\ol{X}|D}$ are clearly functorial for finite co-admissible maps in the sense of \ref{ssec-admissible-coadmissible}.
The subscript $M_{\text{ssup}}$ refers to a slightly different modulus condition on algebraic cycles, stronger then the one considered in \cite{BS}, which we call \textit{strong sup-modulus} condition, following the terminology introduced in \cite{KP} in the case of additive higher Chow groups (see Definition \ref{def:definitionCrXD-both-MC}). Note that for $n=0$, the group $\CH^{d}(\ol{X}|D, 0)_{M_{\text{ssup}}}$ agrees with the Kerz-Saito Chow group of zero cycles with modulus. The appearance of the strong sup-modulus condition, instead of the classical \textit{sum} condition considered elsewhere, is quite natural in our approach, and is related to the fact that the above-mentioned invariance property of relative $K$-theory is built on a certain iterated construction, which forces us to consider ``one face at the time'' of the closed box $\ol{\square}^n = (\mathbb{P}^1)^n$. We expect this new phenomena to play some role in the future development of the theory. 

In this work, we consider the range $(d+n, n)$ for $d=\dim \ol{X}$, i.e. the groups which might be called of zero cycles with coefficients in Milnor $K$-theory. In this sense, the groups $\CH^{d+n}(\ol{X}|D, n)_{M_{\text{ssup}}}$ can be seen as a cycle-theoretic incarnation of the $K$-groups of reciprocity functors of Ivorra-R\"ulling \cite{Ivorra-Rulling} (the $T$-functor, in their terminology) or of the $K$-groups of geometric type considered by Sugiyama \cite{Sugiyama:2017aa}. In fact, as consequence of our Theorem \ref{thm:thm-intro}, we get the following Corollary (see \ref{ssec:Cor-Sugiyama} for a quick review of the notation)
\begin{corintro}[see Corollary \ref{cor:Sugiyama-Cor}] In the notations of Theorem \ref{thm:thm-intro}, assume moreover that $\ol{X}$ is proper over $k$ and that $k$ is perfect. Then there is a canonical homomorphism
    \[{\rm cyc}_{\ol{X}|D, n}^{{\rm geo}} \colon K_{\rm sum}^{\rm geo}(h_0(\ol{X};D), \mathbb{G}_m,   \ldots, \mathbb{G}_m)\tensor\mathbb{Q}\to K_n(\ol{X};D)^{(d+n)}\]
    for each $n\geq 0$. Here $ K_{\rm sum}^{\rm geo}(-)$ denotes Sugiyama's K-group of geometric type of the SC-reciprocity sheaves $h_0(\ol{X};D)$ and $\mathbb{G}_m$ ($n$-times).
    \end{corintro}

In some special cases, e.g. $X=\mathbb{A}^r$ over an algebraically closed field of characteristic zero and $D$ an effective divisor with $\deg(D_{\rm red})<r$, the groups $\CH^{r+n}(\mathbb{A}^r|D,n)$ are known to be zero by results of Krishna-Park \cite[Theorem 1.2 and Theorem 1.3]{KP3}. There are other vanishing results for special fields and for special type of affine varieties, reflecting corresponding vanishing on the $T$-groups. For an example  where the groups  $\CH^{r+n}(\mathbb{A}^r|D,n)$ are \textit{not} zero, see \cite{Crisman2018}. 

Before giving a brief description of the content of the paper, we make a couple of other remarks about our result.
\begin{enumerate}
    \item For our construction, we work with the naive cycle complex rather then with the sheafified version $z^q(-|-\times D, *)_{\rm Nis}$  considered in \cite{IwasaKai} for the definition of relative Chern classes. We expect that the composition of ${\rm cyc}^{d+n}_{\ol{X}|D}$ with the Chern class $c_{d+n, 2d+n}$ agrees, up to a rational factor, with the canonical map $\CH^{d+n}(\ol{X}|D, n)_{M_{\text{ssup}}}\to \CH^{d+n}_{\rm Nis}(\ol{X}|D, n)$. % We need however to assume that the support of $D_{\rm red}$ is a normal crossing divisor in order to perform the required local computations. 
    \item In order to show that the map defined on the set of generators factors through the rational equivalence, we need to consider suitable liftings of classes of $1$-cycles with modulus in relative $K$-theory, but this lift is not canonical enough to be used to produce a cycle class map for $1$-cycles. We explain in \ref{rem:why1-cycles} what is the obstruction for generalizing our method to the construction of cycle classes for higher dimensional cycles. 
    \item We work with a fixed pair $(\ol{X};D)$, and  we do not consider the problem of comparing the pro-groups \[{\rm cyc}^{d+n}_{\ol{X}| mD} \colon\{  \CH^{d+n}(\ol{X}|mD, n)\}_{m, \Q}\to \{K_n(\ol{X}, mD)^{(d+n)}\}_{m, \Q}.\] In the special case $\ol{X} = \Spec{A}$ for $A$ a henselian DVR over $k$ with uniformizer $\pi$ and $D = \Spec{A/\pi}$, a partial result in this direction (but in every range) has been obtained by \cite{IwasaKai} (see Theorem 6.1 in \textit{loc. cit.}).
\end{enumerate}
We finally remark that we expect the maps ${\rm cyc}^{d+n}_{\ol{X}|D}$ to be injective. We verified this for $n=0$ and $k$ algebraically closed when $\ol{X}$ is affine or a quasi-projective surface in \cite[Theorem 1.5 and Theorem 1.7]{BK}. Independently from the results of this paper, an alternative and different construction of a cycle class map for higher zero cycles with modulus to relative K-groups has been recently found by Gupta-Krishna in \cite{Gupta-Krishna}.

\tocless\subsection{} This paper is organized as follows. In Section 1, we review some definitions and elementary results about the cycle complexes with modulus, and we introduce the alternative modulus condition $M_{\rm ssup}$. We refer the reader to \cite{BS} for a more detailed treatment of the subject. In Section 2, we discuss some results on multi-relative $K$-theory, and we explain how to use iteratively the $\mathbb{P}^1$-bundle formula to obtain a delooping of the relative $K$-theory spectrum $\mathbf{K}(\ol{X};D)$. This is what allows us to shift the problem of constructing classes in higher relative $K$-groups to the problem of constructing classes in some co-multi-relative $K_0$-group.

Finally, in Section 3 we present the proof of our main Theorem. Using our new model, the construction of the cycle class map follows the lines of the proof of Levine of the rational isomorphism between Bloch's higher Chow groups and the weight-graded pieces of Quillen's $K$-theory groups (see Section 2 of \cite{LevineBlochrevisited}), but a new argument is necessary to show that the resulting map lifts from $K_n(\ol{X})$ to the relative group $K_n(\ol{X};D)$. This is where the modulus condition on cycles (and in particular the choice of the $M_{\rm ssup}$-condition) comes into the picture, and a delicate local analysis (explained in \ref{sec:lifting-classes}) is required.
 
\subsection*{Notations} Throughout this note, we fix a  field $k$.  All schemes are assumed to be separated and of finite type over $k$ unless otherwise specified. For two $k$-schemes $X$ and $Y$, we write $X\times Y$ for $X\times_k Y$. If $Z$ is any closed subscheme of a scheme $X$, we denote by $|Z|$ its support, i.e.~the underlying closed subspace of $X$. Similarly, if $E$ is an effective Cartier divisor on $X$, we denote by $|E|$ its support.
%\chapter{Algebraic cycles with moduli conditions}
%\chapter{Additive homotopy theory of schemes}

\section{A recollection on relative cycle complexes}\label{sec:rel-cycle-complexes}
 \subsubsection{} Let $\mathbb{P}^1_k = \Proj{k[Y_0, Y_1]}$ be the projective line over $k$ and let $y = Y_1/Y_0$ be the standard rational coordinate function on it. For $n\in \bN\setminus\{0\}, 1\leq i\leq n$, let $p_i^n\colon (\P^1)^n\to \P^1 $ be the projection onto the $i$-th component. We use on $(\P^1)^n$ the rational coordinate system $(y_1,\ldots, y_n)$, where $y_i = y\circ p_i$. Let \[\square^n = (\P^1_k\setminus\{1\})^n\]
    be the \textit{open ($n$-dimensional) box} and let 
$\iota_{i,\epsilon}^n \colon \square^n\to {\square}^{n+1}$ be the closed embedding  
\[\iota_{i,\epsilon}^n(y_1, \ldots, y_n) = (y_1,\ldots, y_{i-1}, \epsilon, y_i, \ldots, y_n), \text{ for } n\in \bN, 1\leq i\leq n+1, \epsilon \in \{0,\infty\},\]
of the codimension one face given by $y_i = \epsilon$, for $\epsilon \in \{0,\infty\}$. The assignment $n \mapsto \square^n$ defines an extended cocubical object $\square^{\bullet}$ in the category of $k$-schemes, in the sense of \cite[1.5]{LevineSmooth}. We conventionally set $\square^{0} =\Spec{k}$.

A face of the open box $\cubb^n$ is a closed subscheme $F$ defined by equations of the form
\[
y_{i_1}=\epsilon_1,\dots,y_{i_r}=\epsilon_r\;;\;\; \epsilon_j\in \{0,\infty\}.
\]
For a face $F$, we write $\iota_F\colon F \hookrightarrow \cubb^n$ for the inclusion. We write $F^n_{i, \epsilon}$ for the face of $\square^n$ given by $y_i = \epsilon$ for $\epsilon \in \{0, \infty\}$ and $i=1,\ldots, n$ (i.e.~the image of $\iota_{i, \epsilon}^{n-1}$ in $\square^n$). The strict normal crossing divisor $\sum_{i, \epsilon} F^n_{i, \epsilon}$ on $\square^n $ will be denoted $\partial \square^n$ and called the \textit{boundary divisor} of $\square^n$.

If no confusion arises, we also write $F^n_{i, \eta}$ for    the divisor of $(\mathbb{P}^1)^n$ given by $y_i = \eta$ for $\eta \in \{0, 1, \infty\}$,  and put
$F^n=\underset{1\leq i\leq n}{\sum} F^n_i$. Finally, we write $\ol{\square}^n$ for $(\P^1)^n$ and we call it the \textit{closed box}.

\subsection{Moduli conditions on higher cycles}\label{sec:moduli-conditions} 
 Let $\Xb$ be a (reduced) equidimensional scheme over $k$, and let $D$ be an effective Cartier divisor on $\Xb$. Let $X$ be the open complement of $D$ in $\Xb$. Let $n >0$ and suppose that $W$ is an integral closed subscheme of $X\times \square^n$. Let $\overline{W}$ denote the closure of $W$ in $\ol{X}\times \overline{\square}^n$ and $\overline{W}^N$ denote its normalization. Write $\varphi_{\overline{W}}\colon \overline{W}^N \to \Xb\times\overline{\square}^n$ for  the composition of the natural map $\overline{W}^{N}\to \overline{W}$ with the closed immersion $\overline{W}\to\Xb\times\overline{\square}^n$.   We say that
\begin{romanlist}
\item $W$ satisfies the $M_{{\sum}}$ modulus condition (the \textit{sum-modulus condition}) if we have the inequality \[\varphi_{\overline{W}}^{*}(D \times \overline{\square}^n)\leq \varphi_{\overline{W}}^{*}(X\times F^{n}).\]
\item  $W$ satisfies the $M_{\text{ssup}}$ modulus condition (the \textit{strong $\sup$-modulus condition}) if there exists an integer $1\leq i\leq n$ such that we have the inequality \[\varphi_{\overline{W}}^{*}(D \times \overline{\square}^n)\leq \varphi_{\overline{W}}^{*}(X\times F^{n}_{i,1}).\] 
\end{romanlist}
The above definitions are generalizations of \cite{KP}. The sum-modulus condition is by far the most used in literature and it's the point of view adopted in \cite{BS}.  Note that the strong-sup condition is strictly stronger then the sum condition (though there are ``conjectures'' about the resulting cycle complexes to be quasi-isomorphic, see again \cite{KP}).
\begin{rmk}There are other possible \textit{moduli conditions} on cycles that are reminiscent of older stages of the theory. For example, in \cite{KP} (generalizing the original definition of Bloch-Esnault \cite{BEAdditive}) one can find a   \textit{sup-modulus condition}, where the relevant inequality of divisors is checked on the supremum over $i$ of $\varphi_{\overline{W}}^{*}(X\times F^{n}_{i,1})$.  
  \end{rmk}
 \subsubsection{}Let $Y$ be a scheme  equidimensional over $k$, $D$ and $F$ two effective Cartier divisors on $Y$. Assume that $D$ and $F$ have no common components. Let $X$ be the open complement $X= Y- (F+D)$. The following Lemma is taken from \cite{BS} and uses the same argument as \cite[Proposition 2.4]{KP}.
  \begin{lem}\label{Containment_YFD} Let $W$ be an integral closed subscheme of $X$ and let $V\subset W$ be an integral closed subscheme of $W$.   Let $\Wb$ (resp. $\Vb$) be the closure of $W$ (resp. of $V$) in $Y$.  Let $\varphi_W\colon \Wb^N\to Y$ (resp. $\varphi_V\colon \Vb^N\to Y$) be the normalization morphism. Then the inequality 
  $\varphi_W^*(D)\leq \varphi_W^*(F)$ as Cartier divisors on $\Wb^N$ implies the inequality
  $\varphi_V^*(D)\leq \varphi_V^*(F)$ as Cartier divisors on $\Vb^N$.
  %Then if $W$ satisfies the modulus condition \ref{moduluscond}, so does $V$.

     \end{lem}
	 
\begin{df}\label{def:definitionCrXD-both-MC}Let $M\in \{M_{{\sum}},  M_{\text{ssup}}\}$. We write $C^r(\Xb|D, n)_M$ for the set of all integral closed subschemes $V$ of codimension $r$ in $X\times{\square}^n$ satisfying the following conditions:
    \begin{enumerate}\item $V$ has proper intersection with $X\times F$ for all faces $F$ of $\square^n$.
        \item For $n=0$, $C^r(\ol{X}|D, 0) = C^r(\ol{X})_D$ is the set of integral closed subschemes of $\ol{X}$ not intersecting $D$.
        \item For $n>0$, $V$ satisfies the $M$-modulus condition on $\ol{X}\times \ol{\square}^n$. %let $\ol{V}$ be the closure of $V$ in $\ol{X}\times \ol{\square}^n$. If $(D\times \ol{\square}^n)\cap \ol{V}\neq \emptyset$, then $V$ satisfies the $M$-modulus condition. 
        \end{enumerate}
        An element in $C^r(\Xb|D, n)_M$ is called an \textit{admissible cycle with $M$-modulus $D$} of codimension $r$, or simply an admissible cycle with modulus $D$ if $M= M_{\sum}$. 
    \end{df}
\begin{rmk}\label{rmk:Modulus-implies-closed-in-complement}Any of the above modulus conditions $M$ imply that 
$\Vb\cap (D\times \ol{\square}^n)\subset \Xb\times F^{n}$ as closed subsets of $\Xb\times \ol{\square}^n$. Hence
$\Vb \cap (D\times \square^n)=\emptyset$ and $V$ is closed in $\Xb\times \square^n$.
This implies that $C^r(\Xb|D, n)_M$ can be viewed as a subset of the set of all integral closed subschemes $W$ of codimension $r$ on $\Xb\times \square^n$ which intersects properly with $\Xb\times F$ for all faces $F$ of $\square^n$.
%For all faces $F\subset \cubb n$ of $\dim(F)>0$ (resp. $\dim(F)=0$), $V$ intersects properly with $D\times F$ (resp. $V\cap (D\times F)=\emptyset$).
\end{rmk}

Let $V\subset W$ be integral closed subschemes of $\ol{X}\times \square^n$ and assume that $W$ satisfies one of the modulus conditions $M\in \{M_{{\sum}},  M_{\text{ssup}}\}$. Lemma \ref{Containment_YFD} shows that the same is true for $V$. This, together with the good position assumption on admissible cycles, proves  the following Lemma.
\begin{lem}\label{lem:containment-applied} Let $V\in C^r(\Xb|D, n)_M$ and let $F$ be a face of dimension $m$ of $\square^n$. Then the cycle $(\id_X\times \iota_F)^*(V)$ on $X\times F \cong X\times \square^m$ belongs to $C^r(\Xb|D,m)_M$.
\end{lem}
\begin{df}\label{def:higher-Chow-modulus-main-def}We denote by $\underline{z}^r(\Xb|D, n)_M$ the free abelian group on the set $ C^r(\Xb|D, n)_M$. By Lemma \ref{lem:containment-applied}, the cocubical object of schemes $n\mapsto \square^n$ gives rise to an extended cubical object (see again \cite[1.5]{LevineSmooth}) in the category of abelian groups
    \[\mathds{1}^n\mapsto \underline{z}^r(\Xb|D, n)_M, \quad n\geq 0.\]
In particular, the groups $ \underline{z}^r(\Xb|D, n)_M$ define a chain complex with boundary
\[\partial = \sum_{1\leq i\leq n} (-1)^i(\partial_i^\infty - \partial_i^0),\]
 where $\partial_i^\epsilon\colon  \underline{z}^r(\Xb|D, n)_M\to  \underline{z}^r(\Xb|D, n-1)_M$ is the pullback along $(\id_{\ol{X}}\times\iota_{i, \epsilon}^{n-1})^*$ for $\epsilon\in \{0, \infty\}$. 
We call the associated non-degenerate complex (see \cite[1.2]{LevineSmooth}) $z^r(\Xb|D, *)_M$ the \textit{cycle complex of $\Xb$ with ($M$)-modulus $D$}.  Its homology groups are denoted by
\[ \CH^r(\Xb|D, n)_M = \H_n(z^r(\Xb|D, *)_M)\]
and called \textit{higher Chow groups of $\Xb$ with ($M$)-modulus $D$}. Note that we have a natural inclusion of cycle complexes \[z^r(\Xb|D, *)_{M_{\rm ssup}} \subset z^r(\Xb|D, *)_{M_\Sigma},\] and therefore natural homomorphisms
\[ \CH^r(\Xb|D, n)_{M_{\rm ssup}} \to \CH^r(\Xb|D,n)_{M_\Sigma} =  \CH^r(\Xb|D,n).\]
    \end{df}
    \begin{rmk} Every admissible cycle with modulus $V\in C^r(\Xb|D, n)_M$ is closed in $\ol{X}\times \square^n$ as noticed in Remark \ref{rmk:Modulus-implies-closed-in-complement}. In particular, we can naturally view the complex $z^r(\Xb|D,*)$ as a subcomplex of the (non-degenerate) cubical cycle complex $z^r(\ol{X},*)$  of Bloch (see \cite[Section 3]{LevineBlochrevisited} for a proof that the cubical version coincides with the original simplicial version of \cite{Bloch86}). This gives a map
        \[ \CH^r(\Xb|D, n)\to \CH^r(\Xb,n) \]
        from higher Chow groups with modulus (for one of the two $M$-conditions) to  Bloch's higher Chow group. Of course, when $D=\emptyset$, there is no modulus condition to check and our definition recovers the usual cubical higher Chow groups.
        \end{rmk}
\begin{rmk}Higher Chow groups with moduli conditions are a generalization of the \textit{additive higher Chow groups} introduced by Bloch-Esnault \cite{BEAdditive} and subsequently studied by Park \cite{ParkAMJ}, R\"ulling \cite{Rul}, Krishna-Levine \cite{KL} and others. For $\Xb = Y\times \mathbb{A}^1_k$, with $Y$ an integral scheme of finite type over $k$, and $D=m\cdot Y\times \{0\} \subset Y\times \mathbb{A}^1_k$ for some $m>0$, the groups $\CH^r(\Xb|D,n)$ coincide with ${\rm TCH}^r(Y, n+1;m)$.
     
    The Definition proposed above, for the $M_\Sigma$ condition, was initially conceived by Kerz and Saito as generalization to higher cycles of the Chow group of zero cycles with modulus used in \cite{KS1} to study wildly ramified class field theory for varieties over finite fields. We refer the reader to \cite{BS} for a systematic treatment.
    \end{rmk}
 
\subsection{Easy functorialities}\label{ssec-admissible-coadmissible}Although we don't need them in this note, we briefly recall some functorial properties of the relative cycle complexes. Let $(\Xb,D)$ and $(\ol{Y}, E)$ be two pairs consisting of reduced equidimensional schemes $\Xb$ and $\ol{Y}$    over $k$ and effective Cartier divisors $D$ and $E$ respectively on them.  When $\Xb$ and $\Yb$ are smooth over $k$, we call the pairs $(\Xb,D)$ and $(\ol{Y}, E)$ \textit{modulus pairs}. 
Let $f\colon Y\to X$ be a morphism in $\mathbf{Sch}(k)$ and assume that the pullback $f^*(D)$ is defined as effective Cartier divisor on $Y$. We say that $f$ is \textit{admissible} (resp.~\textit{coadmissible}) if there is an inequality $f^*(D)\geq E$ (resp.~if there is an inequality $f^*(D)\leq E$). The proofs of the following Lemmas are straightforward.
\begin{lem}\label{lem:pul-flat}Let $f\colon (\ol{Y}, E) \to (\Xb,D)$ be a flat admissible morphism of pairs. The flat pullback of cycles induces a morphism of complexes
    \[f^*\colon  z^r(X| D, *)_M\to z^r(Y|E, *)_M\]
    compatible with composition in the sense that $f^*g^* = (g\circ f)^*$ for  composable admissible flat morphisms  $f$ and $g$. 
\end{lem}
 
\begin{lem}\label{lem:push-proper}Let $f\colon (\ol{Y}, E) \to (\Xb,D)$ be a proper coadmissible morphism of pairs.
 Then there is a well defined push-forward map of cycles
\[f_{*} \colon z^{r+\dim X- \dim Y}(Y| E, *)_M	\to z^r(\Xb|D, *)_M.\]
 
\end{lem}

\section{Relative $K$-theory}\label{sec:cycle-class-zero-cycles}
%Although the sum-modulus condition is by far the most successful in the current development of the theory, there are indeed situations where the strong-sup conditions seems to be better behaved. We present in this Section a construction of a cycle class map from ``higher'' $0$-cycles with modulus to relative $K$-groups.
\subsection{Generalities on relative and multirelative $K$-theory}  
Let $Y$ be a Noetherian separated scheme. Write $\mathbf{K}^{TT}(Y) = \mathbf{K}(Y)$ for the $K$-theory ($\Omega$)-spectrum of Thomason-Trobaugh on the Waldhausen category of strict perfect complexes on $Y$. For a  closed subscheme $j_Z\colon Z\hookrightarrow Y$, the spectrum of algebraic $K$-theory of $Y$ relative to $Z$ is the homotopy fiber of the morphism of spectra $\mathbf{K}(Y) \xrightarrow{j_Z^*} \mathbf{K}(Z)$,
\[\mathbf{K}^{TT}(Y;Z)  = \hofib(\mathbf{K}(Y) \xrightarrow{j_Z^*} \mathbf{K}(Z).)\]
Its homotopy groups, $\pi_*(\mathbf{K}^{TT}(Y;Z)) = K_*(Y;Z)$, are called the \textit{$K$-theory groups of $Y$ relative to $Z$} or simply groups of \textit{relative $K$-theory}. By construction, there is an exact sequence of homotopy groups
\begin{equation}\label{eq:TTrelK} \ldots \to K_{*+1}^{TT}(Z)\to K_{*}^{TT}(Y;Z)\to K_*^{TT}(Y)\to K_*^{TT}(Z)\to K_{*-1}^{TT}(Y;Z)\to \ldots \end{equation}
where $K^{TT}_*(Y)$ and $K_*^{TT}(Z)$ denote, respectively, the Thomason-Trobaugh $K$-theory groups of $Y$ and $Z$. When $Y$ is equipped with an ample family of line bundles as in \cite[Definition 2.1.1]{TT}, the choice of the Waldhausen category is not critical, as remarked by \cite[3.4]{TT}. To simplify the notation, we will drop the superscript $TT$ from the $K$-theory spectra.

\subsubsection{}\label{sec:cube-of-schemes-for-iterated-constructions} Let $T$ be another closed subscheme of $Y$. We denote by $\mathbf{K}^{|T|}(Y)$ or by $\mathbf{K}(Y \text{ on }T)$ the $K$-theory spectrum of the cosimplicial biWaldhausen category of perfect complexes on $Y$ that are acyclic on the open complement $Y\setminus T$. As customary, we call it the \textit{$K$-theory spectrum of $Y$ with support on $T$}. When $U$ is itself quasi compact, Thomason's (proto)-localization theorem \cite[5.1]{TT} gives a homotopy fiber sequence
\[ \mathbf{K}^{|T|}(Y) \to \mathbf{K}(Y) \xrightarrow{\iota_U^*} \mathbf{K}(U) \]
apart possibly from the failure of surjectivity of the map $K_0(Y)\to K_0(U)$. Similarly, for $\cF$ a family of supports on $Y$ we denote by $\mathbf{K}^{\cF}(Y)$ the corresponding $K$-theory spectrum.  In the relative setting, we get the analogue fibration sequence
\[\mathbf{K}^{|T|}(Y;Z) \to \mathbf{K}(Y;Z) \to \mathbf{K}(U;Z\cap U)\]
where the term $\mathbf{K}^{|T|}(Y;Z)$ is defined as the homotopy fiber of the induced map $\mathbf{K}^{|T|}(Y) \to \mathbf{K}^{|T\cap Z|}(Z)$.

\begin{df}\label{def:def-cube-for-iterated-fibers}Let $I$ be a finite set and let $\cP (I)$ be the set of subsets of $I$, seen as  category with morphisms given by inclusions. We call \textit{$I$-cube} a functor $\cX\colon \cP(I)\to \cC$ from $\cP(I)$ to some category $\cC$. An $n$-cube is an $I$-cube with respect to a  set $I$ that has cardinality $n$. Given a subset $J$ of $I$, the inclusion $\cP(J)\to \cP(I)$ defines a $J$-subcube of an $I$-cube $\cX$.
    \end{df}
\subsubsection{} Let $Y$ be a Noetherian separated scheme. Let $Y_1,\ldots, Y_n$ be a set closed subschemes.  We define an $n$-cube of schemes as follows. For every $I\subset \{1, \ldots, n\}$, let $Y_I$ be the subscheme $Y_I = \bigcap_{i\notin I} Y_i$. For $I = \{1, \ldots, n\}$, we conventionally set $Y_{\emptyset} = Y$, and we write $Y_i = Y_{\{1,\ldots, i-1, i+1, \ldots, n\}}$ for short. If $\varphi_{I,J}\colon I\subset J$, there is a corresponding closed embedding of subschemes of $Y$
    \[ \varphi_{I, J} \colon  Y_I =  \bigcap_{i\notin I} Y_i \to \bigcap_{j\notin J} Y_j = Y_J.\]
Pulling back along $\varphi_{I, J}$ defines a (contravariant) $n$-cube of spectra $I\mapsto\mathbf{K}(Y_I)$. 
\begin{df}\label{defn:iterated-fiber} The \textit{total (homotopy) fiber} of $\mathbf{K}(Y_\bullet)$ is, by definition 
     \[ Fib(\mathbf{K}(Y_\bullet)) = \hofib(\mathbf{K}(Y_{{\emptyset}}) \to \holim_{I \neq \emptyset} \mathbf{K}(Y_I)  ). \]
     The \textit{iterated homotopy fiber} of the $K$-theory spectra of $Y_\bullet$ is the $\Omega$-spectrum inductively defined by
     \[\mathbf{K}(Y;Y_1,\ldots, Y_n) = \hofib(\mathbf{K}(Y;Y_1, \ldots, Y_{n-1})\to \mathbf{K}(Y_n; Y_1\cap Y_n,\ldots, Y_{n-1}\cap Y_n)) \]
    \end{df}
  Assume that each intersection scheme $Y_I$ is provided with a family of supports $\cF(Y_I)$ such that, for every $i\in I$ and every $Z\in \cF(Y_I)$, the intersection of $Z$ with $Y_i$ is contained in $\cF(Y_{I\setminus \{i\}})$. We can repeat the above construction replacing everywhere the $K$-theory spectra with the corresponding spectra with support, obtaining in this way a corresponding spectrum $\mathbf{K}^{\cF}(Y;Y_1,\ldots, Y_n)$.
    \subsubsection{}\label{sec:constr-total-iterated-fib-baby}We present now a dual construction. Let $X$ be a Noetherian separated scheme (admitting an ample family of line bundles, as above) and let $i_Y\colon Y\hookrightarrow X$ be a closed subscheme of $X$. Assume that the morphism $i_Y$ is a regular closed immersion. In particular, $i_Y$ is a \textit{perfect} projective morphism in the sense of \cite[Definition 2.5.2]{TT}, and by \cite[3.16.5]{TT} there is a well defined push-forward map $(i_Y)_*\colon \mathbf{K}^{TT}(Y)\to \mathbf{K}^{TT}(X)$.   We denote by $\mathbf{K}(X/Y)$ the homotopy cofiber of $(i_Y)_*$,
    \[ \mathbf{K}(X/Y) = \hocof (\mathbf{K}(Y)\to \mathbf{K}(X)).
    \]
    \subsubsection{}\label{sec:constr-total-iterated-fib}More generally, suppose that we are given a family of closed subschemes $Y_1, \ldots, Y_n$ of $X$. As in \ref{sec:cube-of-schemes-for-iterated-constructions}, consider the $n$-cocube of schemes (see Definition \ref{def:def-cube-for-iterated-fibers}) $Y_\bullet$, setting $Y_{\{1, \ldots, n\}} = X$. 
    Assume now that, for every $I\subset J \subset \{1, \ldots, n\}$, the morphism $\varphi_{I,J}$ is a regular closed immersion. This gives, again by \cite[3.16.5]{TT}, a well-defined push-forward map between the $K$-theory spectra (so, it is a covariant construction)
    \[ (\varphi_{I,J})_*\colon \mathbf{K}(Y_I) \xrightarrow{} \mathbf{K}(Y_J).\]
    We define in this way an $n$-cube of spectra, $\mathbf{K}^{W}(Y_\bullet)$. The \textit{total (homotopy) cofiber}  of $\mathbf{K}(Y_\bullet)$ is by definition
    \[ Cof(\mathbf{K}(Y_\bullet)) = \hocof(\hocolim_{I \neq \{1, \ldots, n\}} \mathbf{K}(Y_I)\to  \mathbf{K}(Y_{\{1, \ldots, n\}})). \]
    \begin{df}\label{def:iter-cof}We keep the above assumptions on $X$ and $Y_1, \ldots, Y_n$. The \textit{iterated homotopy cofiber} of the $K$-theory spectra of $X$ and $Y_*$ is the $\Omega$-spectrum defined inductively
    \[ \mathbf{K}(X / Y_1,\ldots, Y_n) = \hocof ( \mathbf{K}(Y_{n}/ Y_1 \cap Y_{n}, \ldots, Y_{n-1}\cap Y_n) \to \mathbf{K}(X/ Y_1,\ldots, Y_{n-1})).\]
    \end{df}
  \begin{rmk}The following remark holds for the total homotopy fiber as well. There is unique natural map $\mathbf{K}(X / Y_1,\ldots, Y_n) \to Cof(\mathbf{K}(Y_\bullet))$, that is a homotopy equivalence. The existence of the map and the fact that it is an equivalence is dual to \cite[C.6]{FS} (this is also dual to \cite[Proposition 5.5.4]{munson2015cubical}, that the reader can consult for a detailed proof). For every permutation $\sigma$ of the set $\{Y_1,\ldots, Y_n\}$, we have maps
    \[ \mathbf{K}(X / Y_1,\ldots, Y_n) \to Cof(\mathbf{K}(Y_\bullet)) \leftarrow  \mathbf{K}(X / Y_{\sigma(1)},\ldots, Y_{\sigma(n)})\]
    that are homotopy equivalences. In particular, there is a canonical ``zig-zag'' datum joining $\mathbf{K}(X / Y_1,\ldots, Y_n)$ and $\mathbf{K}(X / Y_{\sigma(1)},\ldots, Y_{\sigma(n)})$, and thus the space of homotopies between different iterated homotopy cofibers is contractible. We will then forget the difference between the choices of order of the set of subschemes $Y_*$. \end{rmk}
    \subsection{A model for  relative $K$-theory} We put ourselves back into the geometric situation. Let $Y$ be a regular $k$-variety and consider the $n$-cube of schemes defined by $Y\times \square^n = Y\times (\mathbb{P}^1\setminus \{1\})^n$. Let $\partial \square^n$ denote the strict normal crossing divisor given by the union of the faces $F^n_{i,\varepsilon}$, for $\varepsilon \in \{0, \infty\}$ and $i\leq n$. Using the homotopy property of $K$-theory of regular schemes, there is a natural homotopy equivalence (see \cite[Theorem 3.1]{LevineBlochrevisited})
    \[\mathbf{K}(Y\times\square^n; Y\times \partial \square^n) \to \Omega^n \mathbf{K}(Y)\]
	(where the left-hand side denotes for short the \textit{iterated} homotopy fiber with respect to the boundary divisor $Y\times \partial \square^n$ of $Y\times \square^n$ defined in Section \ref{sec:rel-cycle-complexes})
    giving the isomorphisms 
    \[ K_0(Y\times\square^n; Y\times \partial \square^n)  \xrightarrow{\simeq} K_n(Y), \]
for every $n$. This construction gives a nice delooping of $K$-theory, and allows us to construct classes in higher $K$-groups by constructing classes in (multi)-relative $K_0$.
For $Y$ not regular, the canonical morphism $\mathbf{K}(\mathbb{A}^1_Y)\xrightarrow{\iota_0^*} \mathbf{K}(Y)$ fails to be a homotopy equivalence, and the construction has to be  modified. We take inspiration from \cite{LevineBlochrevisited} in doing so.

The new ingredient is the following: instead of homotopy invariance, we use the projective bundle formula, available by \cite{TT} for any quasi-compact and quasi-separated scheme. We recall the statement. 

\begin{thm}[see \cite{TT}, Theorem 4.1] Let $X$ be a quasi-compact and quasi-separated scheme. Let $\mathcal{E}$ be an algebraic vector bundle of rank $r$ over $X$ and let $\pi\colon \P\mathcal{E}_X\to X$ be the associated projective space bundle. Then there is a natural homotopy equivalence
	\begin{equation}\label{PBF} \prod^r \mathbf{K}(X)\xrightarrow{\sim}\mathbf{K}(\P\mathcal{E}_X)
		\end{equation}
		given by the formula $(x_0, x_1, \ldots, x_{r-1})\mapsto \sum_{i=0}^{r-1}\pi^*(x_i)\tensor[\cO_{\P\mathcal{E}}(-i)]$.
	\end{thm}
\subsubsection{}\label{SectionComputationi_*} Let $X$ be a scheme of finite type over $k$ (though the reader is free to keep working with $X$ quasi-compact and quasi-separated in this subsection).  By \eqref{PBF}, there is an isomorphism
\begin{equation}K_*(\P^1_X)\cong K_*(X)[\cO]\oplus K_*(X)[\cO(-1)],\end{equation}
where $K_*(X)[\cO]$ and $K_*(X)[\cO(-1)]$ are written with respect to the external product
\[K(X)\wedge K(\P^1_\Z)\to K(\P^1_X)\]
and $[\cO]$ and $[\cO(-1)]$ are elements in $K_0(\P^1_\Z)$. It is convenient for us to change basis for the direct sum decomposition to \[\{[\cO], [\cO]- [\cO(-1)]\},\]
so to get 
\begin{equation}\label{PBFP1}K_*(\P^1_X)\cong K_*(X)[\cO]\oplus K_*(X)([\cO]-[\cO(-1)]).\end{equation}
For $i\in \{0,1,\infty\}$, let $\iota_i$ be the regular embedding
\[\iota_{i}\colon X\times \{i\}\to \P^1_X\]
and let $\pi\colon \P^1_X\to X$ be the projection.
We have the associated pullback morphisms
\[\iota_i^*\colon \mathbf{K}(\P^1_X)\to \mathbf{K}(X)\text{ for }i\in \{0,1,\infty\},\]
\[\pi^*\colon \mathbf{K}(X)\to \mathbf{K}(\P^1_X)\]
and the push-forward morphisms
\[\iota_{i,*}\colon \mathbf{K}(X)\to \mathbf{K}(\P^1_X)\text{ for }i\in \{0,1,\infty\}.\]
Note that since the projection $\pi$ has a section, $\pi^*$ is a split monomorphism, corresponding to the canonical inclusion of the first direct summand of \eqref{PBFP1}. For $\iota = \iota_0$, $\iota_\infty$ or $\iota_1$ and for every $j\in \Z$ we have $\iota^*(\cO_{\P^1}(j)) =\cO$. Thus % Let $U(i) = \P^1_X\setminus X\times\{i\}$ for $i\in \{0,\infty\}$. 
\[\iota^*[\cO_{\P^1}] =[\cO] = 1, \text{ and } \iota^*([\cO]-[\cO(-1)]) = 0\] 
in $K_0(X)$. Hence we see that on the direct sum decomposition of $K_*(\P^1_X)$, the pullback morphisms along the  three rational sections all agree and they correspond to the canonical projection on the first component, splitting the pullback along the projection $\pi^*$. In particular, one has that the map
\begin{equation}\label{eq: difference is zero}\mathbf{K}(\P^1_X)\xrightarrow{\iota_0^*-\iota_\infty^*}\mathbf{K}(X)\end{equation}
is homotopy equivalent to the zero map.

By the Projection Formula \cite[3.17]{TT}, the diagram
\[
 \xymatrix{ \mathbf{K}(X)\wedge \mathbf{K}(X)  \ar[r]^{\pi^*\wedge \iota_*} \ar[rd]^{\tensor} & \mathbf{K}(\P^1_X)\wedge \mathbf{K}(\P^1_X) \ar[r]^{\tensor} & \mathbf{K}(\P^1_X) \\
 & \mathbf{K}(X)\ar[ru]_{\iota_*}
 }
    % \begin{diagram}K(X)\wedge K(X) & \rTo^{\pi^*\wedge \iota_*} & K(\P^1_X)\wedge K(\P^1_X) & \rTo^\tensor& K(\P^1_X)\\
    %     &\rdTo_\tensor && \ruTo_{\iota_*} &\\
    %     && K(X)
    % \end{diagram}
\]
commutes, up to canonically chosen homotopy. % Thus the diagram
%\[\xymatrix{ \mathbf{K}(X)\simeq \mathbf{K}[\cO_X] \ar[r]^-{\iota_*} \ar[d]^{\pi^*} & \mathbf{K}(\P^1_X)\\
%\mathbf{K}(\P^1_X)[\iota_* \cO_X]
%} \]
%commutes as well and
 Thus  we see that the push-forward along the inclusion $\iota$ is a monomorphism on the $K$-groups, split by $\pi_*$, that corresponds to the inclusion on the second direct summand of \eqref{PBFP1}, since $[\iota_*\cO_X] = [\cO]-[\cO(-1)] \text{ in }K_0(\P^1_X)$. In particular, the homotopy cofiber \[\mathbf{K}(\P^1_X/ X\times \{1\}) = \mathbf{K}(\P^1_X/ X\times F^1_1)\] is homotopy equivalent to $\mathbf{K}(X)$ via the projection map $\pi^*$.
For $\varepsilon\in \{0, \infty\}$, consider the homotopy fiber
\[\mathbf{K}(\P^1_X; X\times F^1_\varepsilon / X\times F^1_1)  = \hofib( \mathbf{K}(\P^1_X/ X\times F^1_1) \xrightarrow{\iota_\varepsilon^*} \mathbf{K}(X)) \]
Where $F^1_\epsilon = \{\epsilon\} \subset \mathbb{P}^1_k$. Since $\iota_\varepsilon^*$ is a homotopy equivalence, $\mathbf{K}(\P^1_X; X\times F^1_\varepsilon / X\times F^1_1) $ is contractible and thus the iterated homotopy fiber/cofiber
\begin{equation}\label{eq:case-1-delooping}\mathbf{K}(\P^1_X; X\times F^1_0 , X\times F^1_\infty  / X\times F^1_1) \simeq \hocof( \mathbf{K}(X)\xrightarrow{\iota_{1,*}} \mathbf{K}(X\times \P^1; X\times F^1_0 , X\times F^1_\infty)) \end{equation}
is homotopy equivalent to $\Omega \mathbf{K}(X)$.

\subsubsection{}More generally, consider $X\times (\P^1)^n = X\times \ol{\square}^n$. An iterated application of the projective bundle theorem shows that $\mathbf{K}(X\times (\P^1)^n = X\times \ol{\square}^n)$ decomposes as $2n$-copies of the $K$-theory spectrum of $X$, two copies for each copy of $\P^1$ in the closed box $\ol{\square}^n$. 

Let $Cof(\mathbf{K}(X\times (\ol{\square}^n) / X\times F^n_{\bullet,1})$ be total homotopy cofiber of the $n$-cube of schemes \[\{X\times F^n_{i,1} \hookrightarrow X\times (\P^1)^n\}_{i=1}^n,\] where we recall that $F^n_{i,1}$ denotes the face $y_i=1$ on the $i$-th copy of $\P^1$, with respect to the push-forward along the inclusion of faces with value $1$. It is clear by construction and the computation in \ref{SectionComputationi_*} that  $Cof(\mathbf{K}(X\times (\ol{\square}^n) / X\times F^n_{\bullet,1})$ is homotopy equivalent to $\mathbf{K}(X)$.

 For $\emptyset \neq I\subset\{1, \ldots, n\}\times \{0,\infty\}$, consider the subscheme $X\times \partial\ol{\square}^n_I$ of $X\times \ol{\square}^n$ given by 
    \[X\times \partial\ol{\square}^n_I  = \bigcap_{(k, \varepsilon) \not\in I} X\times F^n_{k,\varepsilon} \hookrightarrow X\times \ol{\square}^n\]
For fixed $I$, consider for every  $k$ with $(k,0), (k, \infty)\in I$, the inclusion of the face $\iota_{k,1}^n\colon X\times F^n_{k,1} \to X\times \partial\ol{\square}^n_I$. This defines another (co)cube of schemes, and a corresponding (co)cube of spectra with maps induced by push-forward
    \[\iota_{1,k,*}^n\colon \mathbf{K}(X\times F^n_{k,1}) \to \mathbf{K}(X\times \partial\ol{\square}^n_I = X\times F_{I'}\times\overset{\overset{k}{\vee}}{\P^1}\times F_{I''}) \]
    for a partition $I= I'\cup I''$ with the obvious convention. We denote by $Cof(\mathbf{K}(X\times \partial \ol{\square}^n_I) / X\times F^n_{\bullet,1})$ its total homotopy cofiber. The following Proposition is now proved by descending induction on $n$, starting from \eqref{eq:case-1-delooping} applied to $\ol{X}\times \ol{\square}^n = (\ol{X}\times \ol{\square}^{n-1})\times \mathbb{P}^1$
    \begin{prop}\label{lem:delooping}   The total homotopy fiber/cofiber
\begin{align*}
\mathbf{K}&(X\times \ol{\square}^n; X\times \partial \ol{\square}^n / X\times F^n) \\
 &= \hofib(Cof(\mathbf{K}(X\times \ol{\square}^n) / X\times F^n_{\bullet,1}) \to \holim_{I\neq \emptyset} Cof(\mathbf{K}(X\times \partial \ol{\square}^n_I) / X\times F^n_{\bullet,1}))\end{align*}
        is homotopy equivalent to the $n$-th loop $\Omega^n \mathbf{K}(X)$.
    \end{prop}
\subsubsection{}\label{sssec:relative-setting}
Let $X$ be as above and let $Y$ be a closed subset of $X$ (if the reader is still considering $X$ quasi-compact and quasi-separated, she might want to assume that the open complement $U= X\setminus Y$ is quasi compact as well). The closed immersion $\iota_Y$ gives a pullback morphism on the $K$-theory spectra, and gives induced pullback morphisms between the cubical objects $X\times \overline{\square}^n$ and $Y\times \overline{\square}^n$. We denote by $\mathbf{K}((X;Y)\times \ol{\square}^n; (X;Y)\times \partial \ol{\square}^n / (X;Y)\times F^n) $ the homotopy fiber
	\[\hofib(\mathbf{K}(X\times \ol{\square}^n; X\times \partial \ol{\square}^n / X\times F^n)  \xrightarrow{\iota_Y^*} \mathbf{K}(Y\times \ol{\square}^n; Y\times \partial \ol{\square}^n / Y\times F^n) ).\]
By Proposition \ref{lem:delooping}, we get a natural homotopy equivalence 
	\begin{equation}\label{eq:delooping-final-formula}\mathbf{K}((X;Y)\times \ol{\square}^n; (X;Y)\times \partial \ol{\square}^n / (X;Y)\times F^n) \xrightarrow{\sim}\Omega^n \mathbf{K}(X;Y) \end{equation}
	for the relative $K$ theory spectrum $\mathbf{K}(X;Y)$.

\subsection{Adams operations on $K$-groups and variants} \label{ssec:Adams}  
We collect some useful facts about Adams operations on relative $K$-theory with support. For their construction (and much more), we refer the reader to \cite[Section 5 and 7]{LevineLambdaOp} or to \cite{LecomteAdams}. 

Suppose that $X$ is a Noetherian  scheme equipped with an ample family of line bundles. Then the $K$-groups of $X$, with rational coefficients, come equipped with group homomorphisms (Adams operations)
\[\psi^k\colon K_n(X)_\Q\to K_n(X)_\Q, \quad k\geq 1\] 
which turn the graded ring $\bigoplus_n K_n(X)_\Q$ into a special $\lambda$-ring. By e.g.~\cite[Proposition 4.1.2]{LecomteAdams}, they are functorial for maps of schemes. 

\subsubsection{}When we take homotopy groups, the identification of ${K}_n(\P^1_X/X\times \{1\})_\Q$ with the summand $ K_n(X)_\Q$ of ${K}_n(\P^1_X)_\Q$ via $\pi^*$ is compatible with Adams operations. In fact, the operations $\psi^k$ are compatible with products by \cite[Theorem 4.2.1.iv]{LecomteAdams}, and they satisfy $\psi^k([\cO]) = [\cO]$ by \cite[Theorem 4.2.1.i]{LecomteAdams}. In particular, the decomposition \eqref{PBFP1} is compatible with the action of $\psi^k$ for every $k\geq 1$. 
	
One can extend the argument to the relative $K$-groups, even with support, thanks to the following result of Levine.
\begin{thm}[Corollary 5.6 \cite{LevineLambdaOp}] \label{thm:Levine-lambdaOps}Let $X$ be a Noetherian scheme over a Noetherian ring $S$, admitting an ample family of line bundles. Let $Y_1,\ldots, Y_n$ be closed subschemes of $X$, and let $W$  be another closed subscheme of $X$. Then there is a special $K_0(S)-\lambda$-algebra structure for the relative $K$-theory with supports (see Definition \ref{defn:iterated-fiber})
	\[K_p^{|W|}(X; Y_1,\ldots, Y_n)\]
	which is natural in the tuple $(X; Y_1, \ldots, Y_n; W)$. In particular, the long exact relativization sequence 
	\[\ldots  \to K_p^{|W|}(X; Y_1,\ldots, Y_n ) \to   K_p^{|W|}(X; Y_1,\ldots, Y_{n-1}) \to K_p^{|W\cap Y_n|}(X\cap Y_n; Y_1\cap Y_n, \ldots Y_{n-1}\cap Y_n) \to \ldots \]
	and the long exact localization sequence 
	\[\ldots  \to K_p^{|W|}(X; Y_1,\ldots, Y_n )\to K_p(X; Y_1,\ldots, Y_n ) \to K_p(X\setminus W; Y_1\setminus W,\ldots, Y_n \setminus W) \to \ldots \] 
	are sequences of special $K_0(S)-\lambda$-algebras.	
	\end{thm} 
If we apply Theorem \ref{thm:Levine-lambdaOps} to the case $X\times \P^1$, $Y_1 = X\times F^1_0$ and $Y_2 = X\times F^2_\infty$ (and $S=k$ for a field $k$, say), we get functorial Adams operations on the multi relative $K$-groups $K_n(\P^1_X; X\times F^1_0, X\times F^1_\infty)_\Q$, compatible with the relativization diagram
\[ \xymatrix{ K_n(\P^1_X; X\times F^1_0, X\times F^1_\infty)_\Q \ar[r]\ar[d] &  K_n(\P^1_X; X\times F^1_0)_\Q  \ar[d] \\
  K_n(\P^1_X; X\times F^1_\infty)_\Q \ar[r] & K_n(\P^1_X)_\Q.
 }   \]
 As remarked in \ref{SectionComputationi_*}, the $\P^1$-bundle formula implies that the pullback $\iota_\epsilon^*$, for $\epsilon\in \{0, \infty\}$, factors through the cofiber $K_n(\P^1_X/X\times F^1_1)_\Q$ and gives an isomorphism $K_n(\P^1_X/X\times F^1_1)_\Q \xrightarrow{\iota_\epsilon^*} K_n(X\times F^1_\epsilon)$. Since $\iota_\epsilon^*$ commutes with $\psi^k$ for every $k\geq 1$, this isomorphism is compatible with the Adams operations on both sides (where the operations on the cofiber groups are simply the operations on $K_n(X)_\Q$ transported to $K_n(\P^1_X)_\Q$ via $\pi^*$), and is functorial in $X$. Putting things together, we finally see that the iterated homotopy fiber/cofiber groups 
 \[ K_n(\P^1_X; X\times F_0^1, X\times F_\infty^1/ X\times F_1^1)_\Q\]
 are naturally equipped with functorial (in $X$) Adams operations $\psi^k$, and that the isomorphism 
 \[K_n(\P^1_X; X\times F_0^1, X\times F_\infty^1/ X\times F_1^1)_\Q \to \pi_n (\Omega \mathbf{K}(X))_\Q = K_{n+1}(X)_\Q\]
 given by \eqref{eq:case-1-delooping} is also compatible with Adams operations on both sides.
\subsubsection{}Clearly, the same holds if we replace $\P^1_X$ with $X\times (\P^1)^n$ and we iterate the argument. In particular, thanks to Levine's result, we get functorial Adams operations on the multi-relative fiber/cofiber groups
\[\psi^k\colon K_n(X\times \ol{\square}^n; X\times \partial \ol{\square}^n / X\times F^n)_\Q \to  K_n(X\times \ol{\square}^n; X\times \partial \ol{\square}^n / X\times F^n)_\Q\]
and that the isomorphism given by Proposition \ref{lem:delooping} is compatible with them. 

If now $Y$ is any closed subset of $X$, we can consider the homotopy groups of the homotopy fiber $\mathbf{K}((X;Y)\times \ol{\square}^n; (X;Y)\times \partial \ol{\square}^n / (X;Y)\times F^n)$ introduced in \ref{sssec:relative-setting}. Another application of Theorem \ref{thm:Levine-lambdaOps} gives Adams operations on the $K$-groups 
\[ K_*((X;Y)\times \ol{\square}^n; (X;Y)\times \partial \ol{\square}^n / (X;Y)\times F^n)_\Q\]
compatible with the relativization sequence for the pair $(X;Y)$ and with the isomorphism given by \eqref{eq:delooping-final-formula}. 

\subsubsection{}We finally introduce some some notations. Fix an integer $k>1$. In the setting of Theorem \ref{thm:Levine-lambdaOps}, we denote by \[ K_p^{|W|}(X; Y_1, \ldots, Y_n)^{(q)} \]
the $k^q$-characteristic subspace of $\psi^k$ acting on $K_p^{|W|}(X; Y_1, \ldots, Y_n)_\Q$. This is the set of $v\in K_p^{|W|}(X; Y_1, \ldots, Y_n)_\Q$ such that, for some $N>0$, we have
\[ (\psi^k - k^q\cdot {\rm id})^N(v)=0.\]
The functoriality of the Adams operations on relative $K$-groups with supports shows that the subspaces 
\[K_p^{|W|}(X; Y_1, \ldots, Y_n)^{(q)} \subset K_p^{|W|}(X; Y_1, \ldots, Y_n)_\Q \]
are functorial in the tuple $(X; Y_1,\ldots, Y_n; W)$ (see the remark in \cite[p.~259]{LevineBlochrevisited}). With this notation, the isomorphism \eqref{eq:delooping-final-formula} restricted to the $q$-characteristic subspace takes (on $\pi_0$) the form
\begin{equation}\label{eq:delooping-Adams} K_0((X;Y)\times \ol{\square}^n; (X;Y)\times \partial \ol{\square}^n / (X;Y)\times F^n)^{(q)} \cong K_n(X;Y)^{(q)}. \end{equation}

\def\cycdn{\mathrm{cyc}^{d+n}}
\section{A cycle class map for ``higher'' $0$-cycles with modulus}Assume that $\ol{X}$ is an integral and regular quasi-projective $k$-variety, and let $D$ be an effective Cartier divisor on it. Assume that the support $|D_{\rm red}|$ of $D$ is a strict normal crossing divisor on $\ol{X}$.  We will make systematic  use of the results of \ref{ssec:Adams} about Adams operations on $K$-groups, as well as some other classical facts from \cite[Section 5]{LevineLambdaOp} and \cite[Section 2]{LevineBlochrevisited}.

\subsubsection{} Write $d=\dim \ol{X}$. Recall from Definition \ref{def:higher-Chow-modulus-main-def} that for every $n\geq 0$, the group $z^{d+n}(\ol{X}|D, n)$ is the free abelian group generated by the set $C^{d+n}(\ol{X}|D,n)$ of closed points $P$ in $\ol{X}\times \square^n$ such that $P\notin D\times \square^n$ and $P\notin \ol{X}\times F^{n}_{i, \varepsilon}$ for $i=1,\ldots, n$ and $\epsilon\in \{0, \infty\}$. Clearly, this set coincides with the set of closed points in $\ol{X}\times \ol{\square}^n$ that are disjoint from $D\times \ol{\square}^n$ and that do not meet any face $\ol{X}\times F^n_{i, \eta}$, for $i=1,\ldots,n$ and $\eta\in \{0, 1, \infty\}$.

\subsubsection{}\label{sssec:fund-class} Take a point $P$ in $C^{d+n}(\ol{X}|D,n)$. Since $\ol{X}$ is regular and quasi-projective, the module $\cO_P$ is quasi-isomorphic in the derived category of $\cO_{\ol{X}\times \square^n}$-modules (and in the category of $\cO_{\ol{X}\times \ol{\square}^n}$-modules) to a bounded complex of vector bundles. In particular, we have an isomorphism
\[\mathbb{Q} = K_0(k(P))_{\Q} = K_0(k(P))^{(0)} \xrightarrow{\cong} K_0^{|P|}(\ol{X}\times{\square}^n)^{(d+n)} = K_0^{|P|}(\ol{X}\times \ol{\square}^n)^{(d+n)}.\]
where the last equality follows from the fact that $P\not\in |\sum_{i=1}^n F^n_{i,1}|$. The image of the class of $1$ along the natural morphism $K_0^{|P|}(\ol{X}\times \ol{\square}^n)^{(d+n)}\to K_0(\ol{X}\times \ol{\square}^n)^{(d+n)}$ defines a class $[\cO_P]$, that we call the \textit{fundamental class} of the point $P$. 
\def\Mssup{M_{\rm ssup}}

Let $K_0^{d+n}(\ol{X}\times \ol{\square}^n)^{(d+n)}$ be the direct limit of the groups $K_0^{|P|}(\ol{X}\times \ol{\square}^n)^{(d+n)}$ as $P$ ranges over the set of finite unions of (closed) points $P$ in $C^{d+n}(\ol{X}|D,n)$. Sending $P$ to its fundamental class gives  a  group homomorphism
\begin{equation}\label{eq:first-def-cyc-K0}z^{d+n}(\ol{X}|D, n)_{\mathbb{Q}} \xrightarrow{\cycdn} K_0^{d+n}(\ol{X}\times \ol{\square}^n)^{(d+n)}, \quad \sum_{j=1}^r a_j[P_j]\mapsto \sum_{j=1}^r a_j[\cO_{P_j}]. \end{equation}
Since any $P$ in $C^{d+n}(\ol{X}|D,n)$ is disjoint from $D\times \ol{\square}^n$ and from the boundary divisor $\ol{X}\times \partial \ol{\square}^n$, we have a natural homotopy equivalence between $\mathbf{K}^{|P|}(\ol{X}\times \ol{\square}^n)$ and the multi-relative $K$-spectrum with support $\mathbf{K}^{|P|}((\ol{X};D)\times \ol{\square}^n ; (\ol{X};D)\times \partial\ol{\square}^n)$. The group homomorphism \eqref{eq:first-def-cyc-K0} gives then a group homomorphism to the multi-relative $K_0$-group (defined again as direct limit over $P \in C^{d+n}(\ol{X}|D,n)$ of the $K$-groups with support)
\[z^{d+n}(\ol{X}|D, n)_{\Q} \xrightarrow{\cycdn} K_0^{d+n}((\ol{X};D)\times \ol{\square}^n ; (\ol{X};D)\times \partial\ol{\square}^n)^{(d+n)}.\]
Composing now with the  map induced on $\pi_0$ by the natural morphism of spectra
\[ \mathbf{K}((\ol{X};D)\times \ol{\square}^n; (\ol{X};D)\times \partial \ol{\square}^n )\to   \mathbf{K}((\ol{X};D)\times \ol{\square}^n; (\ol{X};D)\times \partial \ol{\square}^n / (\ol{X};D)\times F^n),\]
and the natural map 
\[ K_0^{d+n}((\ol{X};D)\times \ol{\square}^n ; (\ol{X};D)\times \partial\ol{\square}^n)^{(d+n)} \to K_0((\ol{X};D)\times \ol{\square}^n ; (\ol{X};D)\times \partial\ol{\square}^n)^{(d+n)}\]
we finally obtain a group homomorphism  
\begin{equation}\label{eq:def-of-cycleclass-map} z^{d+n}(\ol{X}|D, n)_{\Q} \xrightarrow{\cycdn} K_0((\ol{X};D)\times \ol{\square}^n; (\ol{X};D)\times \partial \ol{\square}^n / (\ol{X};D)\times F^n)^{(d+n)} \cong K_n(\ol{X};D)^{(d+n)},\end{equation}
where the last isomorphism follows from \eqref{eq:delooping-final-formula}. We will show that this map factors through the higher Chow group $\CH^{d+n}(\ol{X}|D, n)_{\Q,M_{\rm ssup}}$ defined using the strong sup-condition.

\begin{rmk}If we forget about Adams grading, we can still use the fact that we have a chain of  natural isomorphisms  
	\[ \Z = K_0(k(P)) \cong K_0^{|P|}(\ol{X}\times \ol{\square}^n) \cong K_0^{|P|}( (\ol{X}; D)\times \ol{\square}^n ; (\ol{X};D)\times \partial\ol{\square}^n)\]
	to integrally define a group homomorphism
	\[ z^{d+n}(\ol{X}|D, n) \to K_n(\ol{X}; D). \]
	It is not clear, however, if this map factors through the rational equivalence on cycles for  $n> 0$. This is the case when $n=0$, and the interested reader can check that our cycle class map agrees with the one constructed in \cite[Theorem 12.4]{BK}.
	\end{rmk}
\subsection{Exploiting the modulus condition: classes of curves}\label{sec:exploiting-modulus}
We want to study now how to relate a $1$-cycle with modulus with a suitably defined class in the relative $K$-groups. We first discuss how the good-position conditions allow us to construct classes in the relative $K_0$-groups $K_0(\ol{X}\times \ol{\square}^{n+1}; \ol{X}\times \partial \ol{\square}^{n+1})$. This part of the argument is analogue to \cite[Lemma 2.2]{LevineBlochrevisited}.\bigskip

Given any integral curve $C \subset X\times \square^{n+1}$ that is in good position with respect to every face $X\times F^{n+1}_{i, \varepsilon}$, write $\Cb$ for its closure in $\ol{X} \times \ol{\square}^{n+1}$. Let $\cO_{\Cb}$ be its structure sheaf and write $\iota_{\ol{C}}\colon \Cb\to \ol{X}\times \ol{\square}^{n+1}$ for the closed immersion. Since $\ol{X}$ is regular, the coherent $\cO_{\ol{X}\times \ol{\square}^{n+1}}$ module $\cO_{\Cb}$ is quasi-isomorphic to a bounded complex of vector bundles. Suppose moreover that $\Cb$ is itself regular. As for the case of points, sending $1$ to the class of $\cO_{\Cb}$ gives an isomorphism
\[K_0(k(\Cb))_{\Q} = K_0(k(\Cb))^{(0)} \xrightarrow{\cong} K_0^{|{\ol{C}}|}(\ol{X}\times\ol{\square}^{n+1})^{(d+n)}.\]

Suppose now that $W$ is an arbitrary purely $1$-dimensional cycle in $\ol{X}\times \ol{\square}^{n+1}$. Write $z^0(W)_\Q$ for the $\Q$-vector space on the components of $W$. Assume that $W$ is reduced. Removing the $0$-dimensional subset $W'$ of singular points of $W$ does not change the group $z^0(W)_\Q = z^0(W\setminus W')_\Q$. The regularity of $\ol{X}$ gives then an isomorphism  (see the argument at page 263 of \cite{LevineBlochrevisited} )
\begin{equation}\label{eq:K-with-support-and-group-of-cpt-of-a-cycle}z^0(W)_\Q \xrightarrow{\simeq}  K_0^{|W\setminus W'|}(\ol{X}\times \ol{\square}^{n+1}\setminus W')^{(d+n)}\simeq K_0^{|W|}(\ol{X}\times \ol{\square}^{n+1})^{(d+n)}.\end{equation}
Write $z^{d+n}(\ol{X}\times \ol{\square}^{n+1})^W$ for the subgroup of $z^{d+n}(\ol{X}\times \ol{\square}^{n+1})$ supported on $W$, where $z^{d+n}(\ol{X}\times \ol{\square}^{n+1})$ denotes as customary the group of codimension $d+n$ cycles on $\ol{X}\times \ol{\square}^{n+1}$. The isomorphism \eqref{eq:K-with-support-and-group-of-cpt-of-a-cycle} gives then the map
\[z^{d+n}(\ol{X}\times \ol{\square}^{n+1})^W \xrightarrow{\mathrm{cyc}_W}  K_0^{|W|}(\ol{X}\times \ol{\square}^{n+1})^{(d+n)}.\]
\subsubsection{}\label{sec:construction-lifting-boundary}Let now $F$ be a component of $\partial \ol{\square}^n$ and assume that $F$ intersects each component of $W$ properly. We have a commutative diagram
\[\xymatrixcolsep{5pc}\xymatrix{  z^{d+n}(\ol{X}\times \ol{\square}^{n+1})^W_\Q \ar[r]^{\mathrm{cyc}_W} \ar[d]^{\cdot \ol{X}\times F} &  K_0^{|W|}(\ol{X}\times \ol{\square}^{n+1})^{(d+n)}\ar[d]^{\iota_{\ol{X}\times F}^*} \\
 z^{d+n}(\ol{X}\times F)^{W\cap \ol{X}\times F}_\Q  \ar[r]^{\mathrm{cyc}_{W\cap \ol{X}\times F}} & K_0^{|W\cap \ol{X}\times F|}(\ol{X}\times F)^{(d+n)}. 
}
\]
If $[T]$ in $ z^{d+n}(\ol{X}\times \ol{\square}^{n+1})^W_\Q$ is such that $T\cdot (\ol{X}\times F) = 0$, the commutativity of the above diagram implies that the class ${\iota_{\ol{X}\times F}^*}({\rm cyc}_{W}[T])$ is trivial in $ K_0^{|W\cap \ol{X}\times F|}(\ol{X}\times F)^{(d+n)}$. In particular, the class  ${\rm cyc}_{W}[T]$ lifts to the relative $K_0$-group $K_0^{|W|}(\ol{X}\times \ol{\square}^{n+1}; \ol{X}\times F)^{(d+n)}$. Since the $K_1$ group with support $K_1^{|W\cap \ol{X}\times F|}(\ol{X}\times F)^{(d+n)}$ is equal to zero for weight reasons (see \cite[(2.1), p.~261]{LevineBlochrevisited}), this class is well defined.

\begin{notat}\label{notation-deltaprime} We write $\partial'\ol{\square}^{n+1}$ for the divisor $(\sum_{i=2}^{n+1} F^{n+1}_{i, 0} + F^{n+1}_{i, \infty}) + F^{n+1}_{1, \infty}$ and $F^{n+1}$ (resp. $F^n$) for the divisor $\sum_{i=1}^{n+1}F^{n+1}_{i,1}$ (resp.~for the divisor $\sum_{i=1}^{n}F^{n}_{i,1}$) of $\ol{\square}^{n+1}$ (resp.~of $\ol{\square}^{n}$).
    \end{notat}
\subsubsection{}Suppose that $[T]$ satisfies $T\cdot (\ol{X}\times F^{n+1}_{i, \epsilon}) = 0$ for all $i=2,\ldots, n+1$, $\epsilon \in \{0, \infty\}$ and $T\cdot (\ol{X}\times F^{n+1}_{1, \infty}) = 0$.  We can iterate the argument of \ref{sec:construction-lifting-boundary} to get inductively a well defined class 
\[{\rm cyc}_{W}[T] \in K_0^{|W|}(\ol{X}\times \ol{\square}^{n+1}; \ol{X}\times \partial'\ol{\square}^{n+1})^{(d+n)}.\]
Projecting to the iterated cofiber along the faces $F^{n+1}_{i, 1}$ for $i=1,\ldots, n+1$, gives then a class (that we still denote in the same way)
\[{\rm cyc}_{W}[T] \in K_0^{|W|}(\ol{X}\times \ol{\square}^{n+1}; \ol{X}\times \partial'\ol{\square}^{n+1} / \ol{X}\times F^{n+1})^{(d+n)}\]
and forgetting the support we end up with a class 
\[{\rm cyc}_{d+n}[T] \in K_0(\ol{X}\times \ol{\square}^{n+1}; \ol{X}\times \partial'\ol{\square}^{n+1} / \ol{X}\times F^{n+1})^{(d+n)}.\]
%We will now see how to use the modulus condition to further lift the class with support ${\rm cyc}_{W}[T]$ to a class in the (multi)relative $K_0$-group of $(\ol{X};D)$.
\subsubsection{}\label{sec:notation-classes-to-lift}
Let $N(z^{d+n}(\ol{X}|D, n+1)_{M_{\rm ssup}})$ be the group of admissible cycles in the normalized complex of $\underline{z}^{d+n}(\ol{X}|D,\bullet)_{M_{\rm ssup}}$ (see \cite[p.7]{LevineSmooth} for the definition of normalized subcomplex associated to a cubical object in an abelian category). To simplify this already heavy notation, we suppress the subscript $M_{\rm ssup}$ in what follows.

A cycle in  $N(z^{d+n}(\ol{X}|D, n+1))$ is a $1$-dimensional cycle $Z$ in $X\times \square^n$ such that, for every face $F^{n+1}_{i, \varepsilon}$, for $i=1,\ldots, n+1$, $\epsilon\in\{0, \infty\}$ but with $(i, \epsilon)\neq(1,0)$, it satisfies $Z\cdot F^{n+1}_{i,\varepsilon} = 0$. Moreover, $Z$ is in good position with respect to the remaining face $X\times F^{n+1}_{1, 0}$ and it satisfies the $M_{\rm ssup}$ modulus condition. We can furthermore assume that no component of $Z$ is a vertical coordinate line, i.e., the pullback along a projection $p_j\colon X\times \square^{n+1}\to X\times \square^n$ of a point $P\in X\times \square^n$.

The group of $0$-cycles with modulus $\CH^{d+n}(\ol{X}|D,n)_{\Mssup}$ is then by \cite[Lemma 1.6]{LevineSmooth} the cokernel
\[N(z^{d+n}(\ol{X}|D, n+1))_\Q \xrightarrow{\cdot X\times F^{n+1}_{1,0}} z^{d+n}(\ol{X}|D,n)_\Q\to \CH^{d+n}(\ol{X}|D,n)_{\Q,\Mssup} \to 0. \]
\subsubsection{}\label{sec:notation-classes-to-lift-2} Let $Z$ be a normalized admissible cycle $Z = \sum_{k=1}^r m_k Z_k\in N(z^{d+n}(\ol{X}|D, n+1))$. The closure $\ol{Z}$ of $Z$ in $\ol{X}\times \ol{\square}^{n+1}$ is the closure of its components $Z_1,\ldots, Z_r$. Note that the class 
\[{\rm cyc}_{\ol{Z}}[\ol{Z}_k]\in K_0^{|\ol{Z}|}(\ol{X}\times \cub {n+1})^{(d+n)}\] 
is the image of ${\rm cyc}_{\ol{Z}_k}[\ol{Z}_k] \in K_0^{|\ol{Z}_k|}(\ol{X}\times \cub {n+1})^{(d+n)}$ via the natural map
\[\rho_{\ol{Z}_k, \ol{Z}}\colon  K_0^{|\ol{Z}_k|}(\ol{X}\times \cub {n+1})^{(d+n)} \to K_0^{|\ol{Z}|}(\ol{X}\times \cub {n+1})^{(d+n)}\cong\bigoplus_{k=1}^r K_0^{|\ol{Z}_k|}( \ol{X}\times \cub {n+1})^{(d+n)}  \]
where isomorphism follows from \eqref{eq:K-with-support-and-group-of-cpt-of-a-cycle}. In particular, the class ${\rm cyc}_{\ol{Z}}[\ol{Z}]$ is the sum $\sum_{k=1}^r m_k\rho_{\ol{Z}_k, \ol{Z}} {\rm cyc}_{\ol{Z}_k}[\ol{Z}_k]$. 
By definition of the $M_{\rm ssup}$-modulus condition, each $\ol{Z}_k$ satisfies the condition \ref{sec:moduli-conditions} ii) with respect to some face $j=j(k)\in \{1, \ldots, n+1\}$. Given this, we will need to show the following
\begin{claim}\label{main-claim-vanishing-class} With the above notations, assume that $|D_{\rm red}|$ is a strict normal crossing divisor on $\ol{X}$. Then the image of the class  ${\rm cyc}_{\ol{Z}}[\ol{Z}_k]$ in the cofiber group \[K_0^{|\ol{Z}|}(\ol{X}\times \cub {n+1} / \ol{X}\times F^{n+1}_{j,1})^{(d+n)}\] vanishes along the restriction to $ K_0^{|\ol{Z}|}(D\times \cub {n+1}  / D\times F^{n+1}_{j,1})^{(d+n)}$.
    \end{claim}
Note that the required vanishing can be checked on $K_0$-group with support in $\ol{Z}_k$, and that this immediately implies the same vanishing in the bigger group with support in $\ol{Z}$. This Claim will allow us to lift the image of ${\rm cyc}_{\ol{Z}}[\ol{Z}_k]$ to a class in the relative group 
 \begin{equation}\label{eq:relative-fiber-cofiber-j}K_0^{|\ol{Z}|}((\ol{X};D)\times \cub {n+1} / \ol{X}\times F^{n+1}_{j,1})^{(d+n)}\end{equation}
modulo the image of $K_1^{|\ol{Z}\cap D\times \cub {n+1}|}(D\times \cub {n+1}  / D\times F^{n+1}_{j,1})^{(d+n)}$. In the proof  we will see how the modulus condition on the cycle plays a substantial role.

\begin{rmk}\label{rem:why1-cycles} Since $D$ is not regular, we \emph{cannot} conclude as before that \[K_1^{|\ol{Z}|}(D\times \cub {n+1}  / D\times F^{n+1}_{j,1})^{(d+n)}=0.\] In fact, the vanishing of $K_1^{|W\cap \ol{X}\times F|}(\ol{X}\times F)^{(d+n)}$ in  \ref{sec:construction-lifting-boundary} is a special case of \cite[Claim (2.1)]{LevineBlochrevisited}, that uses the regularity of $\ol{X}$ and of the face $F$ in an essential way. Since the class we are after is necessary only to produce relations in the relative $K_0$, we will not worry about the problem of the choice of the lifting. Of course, this would be the first problem to solve in order to construct a cycle class map for $1$-cycles with modulus. See Remark \ref{rmk:charp} below for a comment in positive characteristic.
    \end{rmk}
We postpone the proof of Claim \ref{main-claim-vanishing-class} to  Section \ref{sec:lifting-classes}. Write $p_{\ol{X}}^j$ for the natural map  
 \[p_{\ol{X}}^j\colon K_0^{|\ol{Z}|}(\ol{X}\times \cub {n+1})^{(d+n)}\to K_0^{|\ol{Z}|}(\ol{X}\times \cub {n+1} / \ol{X}\times F^{n+1}_{j,1})^{(d+n)} \]
 and    $p_{\ol{X}}^{j,i}$  for the map 
 % \[ p_{\ol{X}}^{i,j} \colon K_0^{|\ol{Z}|}(\ol{X}\times \cub {n+1} / \ol{X}\times F^{n+1}_{i,1})^{(d+n)} \to  K_0^{|\ol{Z}|}(\ol{X}\times \cub {n+1} / \ol{X}\times (F^{n+1}_{j,1}, F^{n+1}_{i,1} ))^{(d+n)}\]
 \[p_{\ol{X}}^{j,i}\colon  K_0^{|\ol{Z}|}(\ol{X}\times \cub {n+1} / \ol{X}\times F^{n+1}_{j,1})^{(d+n)} \to K_0^{|\ol{Z}|}(\ol{X}\times \cub {n+1} / \ol{X}\times (F^{n+1}_{j,1}, F^{n+1}_{i,1} ))^{(d+n)}.\]
 Denote similarly by   $p_{D}^{j,i}$ the corresponding map  in the $K$-groups of $D$. We can trace the image of the class ${\rm cyc}_{\ol{Z}}[\ol{Z}_k]$ in the iterated cofiber group 
  \[K_0^{|\ol{Z}|}(\ol{X}\times \cub {n+1} / \ol{X}\times F^{n+1})^{(d+n)}\]
 and its restriction to the corresponding group for $D$ as follows. For $i\neq j$ we look at the commutative diagram
 \[\xymatrix{ K_0^{|\ol{Z}|}(\ol{X}\times \cub {n+1} / \ol{X}\times F^{n+1}_{j,1})^{(d+n)} \ar[r]^{{\rm res}_j} \ar[d]^{p_{\ol{X}}^{j,i}}  &  K_0^{|\ol{Z}\cap D\times \ol{\square}^{n+1}|}(D\times \cub {n+1} / D\times F^{n+1}_{j,1})^{(d+n)}\ar[d]^{p_{D}^{j,i}} \\
  K_0^{|\ol{Z}|}(\ol{X}\times \cub {n+1} / \ol{X}\times (F^{n+1}_{j,1}, F^{n+1}_{i,1} ))^{(d+n)} \ar[r]^-{{\rm res}_{j,i}} & K_0^{|\ol{Z}\cap D\times \ol{\square}^{n+1}|}(D\times \cub {n+1} / D\times (F^{n+1}_{j,1}, F^{n+1}_{i,1} ))^{(d+n)},\\
  %K_0^{|\ol{Z}|}(\ol{X}\times \cub {n+1} / \ol{X}\times F^{n+1}_{i,1})^{(d+n)} \ar[r]^{{\rm res}_i} \ar[u]^{p_{\ol{X}}^{i,j}}  &  K_0^{|\ol{Z}\cap D\times \ol{\square}^{n+1}|}(D\times \cub {n+1} / D\times F^{n+1}_{i,1})^{(d+n)}\ar[u]^{p_{D}^{i,j}}
 }
 \]
where ${\rm res}_j$ and ${\rm res}_{j,i}$ are the obvious restriction maps (induced by $\iota^*_{D\times \ol{\square}^{n+1}}$). Suppose that the component $Z_k$ of the cycle $Z$ satisfies the modulus condition with respect to the face $j$, as stated in Claim  \ref{main-claim-vanishing-class}. Then, the restriction ${\rm res}_j( p^j_{\ol{X}}{\rm cyc}_{\ol{Z}}[\ol{Z}_k])$ vanishes. Since the above diagram commutes, we see that the same thing holds for the restriction
\[ {\rm res}_{j,i} ( p_{\ol{X}}^{j,i} ( p_{\ol{X}}^j  ({\rm cyc}_{\ol{Z}}[\ol{Z}_k]  ))) =0 \quad \text{ in }  K_0^{|\ol{Z}\cap D\times \ol{\square}^{n+1}|}(D\times \cub {n+1} / D\times (F^{n+1}_{j,1}, F^{n+1}_{i,1} ))^{(d+n)}. \]
  Repeating this process for every face $F^{n+1}_{i,1}$, $i=1,\ldots, n+1$, we obtain that the projection of the class ${\rm cyc}_{\ol{Z}}[\ol{Z}_k]$ in the iterated cofiber group
  \[\bigoplus_{k=1}^r K_0^{|\ol{Z}_k|}( \ol{X}\times \cub {n+1})^{(d+n)} \cong K_0^{|\ol{Z}|}(\ol{X}\times \cub{n+1})^{(d+n)} \to K_0^{|\ol{Z}|}(\ol{X}\times \cub{n+1} / \ol{X}\times F^{n+1})^{(d+n)}  \]
dies in the group  $K_0^{|\ol{Z}|}(D\times \cub{n+1} / D\times F^{n+1})^{(d+n)}$ (which now does not depend on the face $j=j(k)$ where the modulus condition for $Z_k$ is fulfilled). Since the class ${\rm cyc}_{\ol{Z}}[\ol{Z}]$ is the sum of the classes $\sum_{k=1}^{r} m_k  {\rm cyc}_{\ol{Z}} [\ol{Z_k}]$, each $Z_k$ satisfies the modulus condition with respect to some face $F^{n+1}_{j(k),1}$,  and the restriction map
\[  {\rm res}^{n+1}_{\ol{Z}}\colon K_0^{|\ol{Z}|}(\ol{X}\times \cub{n+1} / \ol{X}\times F^{n+1})^{(d+n)} \to K_0^{|\ol{Z}\cap D\times \cub {n+1}|}( D\times \cub{n+1} / D\times F^{n+1})^{(d+n)}\]
 is obviously a group homomorphism, we finally see that we have 
 \[{\rm res}^{n+1}_{\ol{Z}} ( {\rm cyc}_{\ol{Z}}[\ol{Z}]) = \sum_{k=1}^r m_k {\rm res}^{n+1}_{\ol{Z}}( {\rm cyc}_{\ol{Z}} [\ol{Z_k}]) =0.\]
In particular, from the long exact relativization sequence, we obtain a class $\beta(\ol{Z})$ in the iterated relative cofiber
 \[ \beta(\ol{Z}) \in K_0^{|\ol{Z}|}((\ol{X};D)\times \cub {n+1} / (\ol{X};D)\times F^{n+1})^{(d+n)}, \]
 well defined up to a class in the image of the natural map
 \[K_1^{|\ol{Z}\cap D\times \ol{\square}^{n+1}|}(D\times \cub {n+1} / D\times F^{n+1})^{(d+n)} \to K_0^{|\ol{Z}|}((\ol{X};D)\times \cub {n+1} / (\ol{X};D)\times F^{n+1})^{(d+n)}, \]
and lifting the fundamental class  ${\rm cyc}_{\ol{Z}}[\ol{Z}]$ of \eqref{eq:K-with-support-and-group-of-cpt-of-a-cycle} along the map 
\[ K_0^{|\ol{Z}|}((\ol{X};D)\times \cub {n+1} / (\ol{X};D)\times F^{n+1})^{(d+n)} \to K_0^{|\ol{Z}|}(\ol{X}\times \cub{n+1} / \ol{X}\times F^{n+1})^{(d+n)}.\]

\subsubsection{}\label{sec:the-class-does-not-change-much} By assumption, the cycle $Z$ satisfies $Z\cdot ({X}\times F^{n+1}_{i, \epsilon}) = 0$ for all $i=2,\ldots, n+1$, $\epsilon \in \{0, \infty\}$ and ${Z}\cdot ({X}\times F^{n+1}_{1, \infty}) = 0$. By Remark \ref{rmk:Modulus-implies-closed-in-complement}, $Z$ is already closed in $\ol{X}\times \square^{n+1}$, so that any extra point of intersection of $\ol{Z}$ with a face $\ol{X}\times \ol{F}^{n+1}_{i, \epsilon} \subset \ol{X}\times \ol{\square}^{n+1}$ is supported on some intersection $\ol{X}\times \ol{F}^{n+1}_{i, \epsilon}\cap \ol{X}\times F^{n+1}_{k, 1}$ for some $k\in \{1,\ldots, n+1\}$, $k\neq i$ (we introduce the overline notation $\ol{F}^{n+1}_{i, \epsilon}$ for sake of clarity).

In particular, the class $\iota^*_{\ol{X}\times F^{n+1}_{i, \epsilon}}({\rm cyc}[\ol{Z}])$ is trivial in the quotient
\begin{equation}\label{eq:quotient-cofiber-regular} K_0^{|\ol{Z}\cap\ol{X}\times F^{n+1}_{i, \epsilon} |}(\ol{X}\times F^{n+1}_{i, \epsilon}/ \ol{X}\times F^{n+1}_{k,1} \cap \ol{X}\times F^{n+1}_{i,\epsilon} )^{(d+n)}.\end{equation}
\begin{rmk}The scheme $\ol{X}\times F^{n+1}_{i, \epsilon}$ is regular, and the $K_0$-group with support on $\ol{Z}\cap\ol{X}\times F^{n+1}_{i, \epsilon}$ depends only on the set of points $|\ol{Z}\cap\ol{X}\times F^{n+1}_{i, \epsilon}|$. In particular, we have  isomorphisms 
    \[ K_0^{|\ol{Z}\cap\ol{X}\times F^{n+1}_{i, \epsilon}|}(\ol{X}\times F^{n+1}_{i, \epsilon}))^{(d+n)} \cong K_0(k(\ol{Z}\cap\ol{X}\times F^{n+1}_{i, \epsilon}))^{(0)} \cong  K_0(k(\ol{Z}\cap\ol{X}\times F^{n+1}_{i, \epsilon}\cap \ol{X}\times F^{n+1}_{k, 1}))^{(0)}\]
    showing immediately that $\iota^*_{\ol{X}\times F^{n+1}_{i, \epsilon}}({\rm cyc}[\ol{Z}])$ dies in the quotient \eqref{eq:quotient-cofiber-regular}. This argument fails  without the regularity assumption on $\ol{X}$, since we can't identify the $K_0$ with support with its support. This is precisely the problem of lifting classes to the group \eqref{eq:relative-fiber-cofiber-j} that we discussed before. %, and that we will solve in Section \ref{sec:lifting-classes} using the modulus condition.
    \end{rmk}
The vanishing of $\iota^*_{\ol{X}\times F^{n+1}_{i, \epsilon}}({\rm cyc}[\ol{Z}])$ in the quotient group \eqref{eq:quotient-cofiber-regular} allows us to lift it to a class in the relative group
\[ {\rm cyc}[\ol{Z}]\in K_0^{|\ol{Z}|}(\ol{X}\times \ol{\square}^{n+1}; \ol{X}\times F^{n+1}_{i,\epsilon}/ \ol{X}\times F^{n+1}_{k,1})^{(d+n)}\]
well-defined, since the group $K_1^{|\ol{Z}\cap\ol{X}\times F^{n+1}_{i, \epsilon} |}(\ol{X}\times F^{n+1}_{i, \epsilon}/ \ol{X}\times F^{n+1}_{k,1} \cap \ol{X}\times F^{n+1}_{i,\epsilon} )^{(d+n)}$ is trivial (in fact, the group $K_1^{|\ol{Z}\cap\ol{X}\times F^{n+1}_{i, \epsilon} |}(\ol{X}\times F^{n+1}_{i, \epsilon})^{(d+n)}$ is trivial again by weight reasons, using again \cite[(2.1), p.~261]{LevineBlochrevisited}, and so is, \textit{a fortiori}, the cofiber group). We repeat the argument for every $i=2,\ldots, n+1$ and $\epsilon\in \{0, \infty\}$ and one more time for $i=1$, $\epsilon = \infty$ to get a class in the iterated cofiber-fiber
\[{\rm cyc}[\ol{Z}]\in K_0^{|\ol{Z}|}(\ol{X}\times \ol{\square}^{n+1}; \ol{X}\times \partial'\ol{\square}^{n+1}/ \ol{X}\times F^{n+1})^{(d+n)} \]
where $\partial'\ol{\square}^{n+1}$ and $F^{n+1}_{k,1}$ are defined as in Notation \ref{notation-deltaprime}. 
% denotes as above the divisor $(\sum_{i=2}^{n+1} F^{n+1}_{i, 0} + F^{n+1}_{i, \infty}) + F^{n+1}_{1, \infty}$ and $F^{n+1}$ is the sum of all the extra faces $F^{n+1}_{k,1}$. 

The lifting property that we just showed for ${\rm cyc}[\ol{Z}]$ implies in fact the same 
property for the chosen ``relative'' lift $\beta(\ol{Z})$: this is in fact obvious once one accepts our claim about the vanishing of the restriction of ${\rm cyc}[\ol{Z}]$ to the cofiber group $K_0^{|\ol{Z}|}(D\times \ol{\square}^{n+1}/ D\times F^{n+1} )^{(d+n)}$.
 
In particular, we finally obtain a class~---~that we keep denoting $\beta(\ol{Z})$~---~in the following multi-relative $K$-group:
\begin{equation}\label{eq:fundamental-class-beta-with-support} \beta(\ol{Z})\in K_0^{|\ol{Z}|}((\ol{X};D)\times \ol{\square}^{n+1}; (\ol{X};D)\times \partial'\ol{\square}^{n+1}/(\ol{X};D)\times F^{n+1})^{(d+n)},\end{equation}
which is well defined up to a relative $K_1$-class, supported on $|\ol{Z}\cap D\times \ol{\square}^{n+1}|$. 
Forgetting the support gives a class
\[\beta_{d+n}[\ol{Z}] \in K_0((\ol{X};D)\times \ol{\square}^{n+1}; (\ol{X};D)\times \partial'\ol{\square}^{n+1}/(\ol{X};D)\times F^{n+1})^{(d+n)}\]
well defined up to elements in the image of $K_1(D\times \ol{\square}^{n+1}; D\times \partial'\ol{\square}^{n+1}/D\times F^{n+1})^{(d+n)}$. 

We can actually give a more precise statement. Let $G_0^{d+n}(\ol{X}|D, n+1)$  be the direct limit 
        \[ G_0^{d+n}(\ol{X}|D, n+1) = \varcolim_{Z\in N(z^{d+n}(\ol{X}|D, n+1))_\Q} K_0^{|\ol{Z}|}((\ol{X};D)\times \ol{\square}^{n+1}; (\ol{X};D)\times \partial'\ol{\square}^{n+1}/(\ol{X};D)\times F^{n+1})^{(d+n)}.
        \] If we write $\Phi_{Z}$ for the composite map 
\[ \xymatrix{K_1^{|\ol{Z}\cap D\times \ol{\square}^{n+1}|}(D\times \ol{\square}^{n+1}; D\times \partial'\ol{\square}^{n+1}/D\times F^{n+1})^{(d+n)}\ar[d]   \\ K_0^{|\ol{Z}|}((\ol{X};D)\times \ol{\square}^{n+1}; (\ol{X};D)\times \partial'\ol{\square}^{n+1}/(\ol{X};D)\times F^{n+1})^{(d+n)} \ar[d] \\
G_0^{d+n}(\ol{X}|D, n+1)   \ar[d]\\
 K_0((\ol{X};D)\times \ol{\square}^{n+1}; (\ol{X};D)\times \partial'\ol{\square}^{n+1}/(\ol{X};D)\times F^{n+1})^{(d+n)}  
}\]
we see that the class $\beta_{d+n}[\ol{Z}]$ is well defined up to the image of $\Phi_{Z}$. %The advantage of this description is explained by the following
%\begin{lem} Every class in the image of $\Phi_{Z}$ goes to zero in $K_0((\ol{X};D)\times \ol{\square}^n; (\ol{X};D)\times \partial \ol{\square}^n / (\ol{X};D)\times F^n)^{(d+n)}$ under the restriction map $\iota^*_{\ol{X}\times F^{n+1}_{1,0}}$.
    % \begin{proof} This is another feature of the modulus condition: the support of $\iota^*_{\ol{X}\times F^{n+1}_{1,0}}(\ol{Z}\cap D\times \ol{\square}^{n+1})$, if it is not already empty, has to be contained in the intersection  $\ol{X}\times F^{n+1}_{1,0}\cap \ol{X}\times F^{n+1}_{j,1}$. It is clear then that such class has to vanish in the iterated cofiber $K_0((\ol{X}, D)\times \ol{\square}^{n}/(\ol{X};D)\times F^n)$.
    %     \end{proof}
%\end{lem}
Write $G_1^{d+n}(\ol{X}|D,n+1)$ for the limit
\[ G_1^{d+n}(\ol{X}|D,n+1) = \varcolim_{Z\in N(z^{d+n}(\ol{X}|D, n+1))_\Q}K_1^{|\ol{Z}\cap D\times \ol{\square}^{n+1}|}(D\times \ol{\square}^{n+1}; D\times \partial'\ol{\square}^{n+1}/D\times F^{n+1})^{(d+n)} \]
and let $\Phi\colon G_1^{d+n}(\ol{X}|D,n+1) \to  K_0((\ol{X};D)\times \ol{\square}^{n+1}; (\ol{X};D)\times \partial'\ol{\square}^{n+1}/(\ol{X};D)\times F^{n+1})^{(d+n)}$ be the induced map. Note that it clearly factors through $G_0^{d+n}(\ol{X}|D, n+1)$ (and we still denote by $\Phi$ the map obtained in this way). We can at this point define a map 
\[ N(z^{d+n}(\ol{X}|D, n+1))_\Q  \xrightarrow{\beta_{d+n}} K_0^{d+n}(\ol{X}|D, n+1) \] 
% \widetilde{K_0}((\ol{X};D)\times \ol{\square}^{n+1}; (\ol{X};D)\times \partial'\ol{\square}^{n+1}/(\ol{X};D)\times F^{n+1})^{(d+n)} 
where  $K_0^{d+n}(\ol{X}|D, n+1)$ the quotient of $G_0^{d+n}(\ol{X}|D, n+1)$ by the image of $\Phi$. % $\widetilde{K_0}$ denotes the quotient of the multi-relative $K_0$ group by the image of $\Phi$.
 Using this construction, we can finally prove the following
\begin{thm}\label{thm:main-theorem-cycle-class}Let $\ol{X}$ be a smooth, quasi-projective $k$-variety of dimension $d$ and let $D$ be an effective Cartier divisor on $\ol{X}$ such that $|D_{\rm red}|$ is a strict normal crossing divisor. Then the cycle class map
    \[z^{d+n}(\ol{X}|D, n)_{\Q} \xrightarrow{\cycdn} K_0((\ol{X};D)\times \ol{\square}^n; (\ol{X};D)\times \partial \ol{\square}^n / (\ol{X};D)\times F^n)^{(d+n)} \cong K_n(\ol{X};D)^{(d+n)}\]
    factors through $\CH^{d+n}(\ol{X}|D,n)_{\Q}$.  Here the last isomorphism follows from \eqref{eq:delooping-Adams}.
    \begin{proof} 
We have to show that the composition  $\cycdn \circ(- \cdot X\times F^{n+1}_{1,0})$ is zero.  Let $Z\in  N(z^{d+n}(\ol{X}|D, n+1))_\Q$ and choose a class $\beta(\ol{Z})$ as in \eqref{eq:fundamental-class-beta-with-support}, lifting the canonical class ${\rm cyc}[\ol{Z}]$. Let $W = \ol{Z}\cap \ol{X}\times \ol{F}^{n+1}_{1,0}$ and $W' =  Z\cap \ol{X}\times F^{n+1}_{1,0}$ (both with reduced scheme structure). Clearly, the finite set of points $W\setminus W'$ is contained in the union of the faces $F^{n+1}_{k,1}$ of $\ol{\square}^{n+1}$. We have
        \[ \xymatrix{K_0^{|\ol{Z}|}(\ol{X}\times \ol{\square}^{n+1})^{(d+n)} \ar[r] \ar[d]^{\iota^*_{\ol{X}\times \ol{F}^{n+1}_{1,0}}}  & K_0^{|\ol{Z}|}(\ol{X}\times \ol{\square}^{n+1} / \ol{X}\times F^{n+1})^{(d+n)} \ar[d]^{\iota^*_{\ol{X}\times \ol{F}^{n+1}_{1,0}}}  & \\
        K_0^{|W|}(\ol{X}\times \ol{F}^{n+1}_{1,0})^{(d+n)} \ar[r] & K_0^{|W|}(\ol{X}\times \ol{F}^{n+1}_{1,0} / \ol{X}\times F^{n})^{(d+n)} \ar[r]^{\cong} &  K_0^{|W'|}(\ol{X}\times F^{n+1}_{1,0} / \ol{X}\times F^{n})^{(d+n)}
        }
        \]
        Where the last isomorphism follows from the same argument used in \ref{sec:the-class-does-not-change-much}.
In particular, the class of $\iota^*_{\ol{X}\times \ol{F}^{n+1}_{1,0}}({\rm cyc}[\ol{Z}])$ in the cofiber group $K_0^{|W|}(\ol{X}\times \ol{F}^{n+1}_{1,0} / \ol{X}\times F^{n})^{(d+n)}$ agrees with the class 
\[\iota^*_{\ol{X}\times F^{n+1}} ({\rm cyc}[Z])\in K_0^{|W'|}(\ol{X}\times \ol{\square}^{n+1} /\ol{X}\times F^{n})^{(d+n)},\]
 so that we don't see the ``extra intersection points'' given by the closure of $Z$ in $\ol{X}\times \ol{\square}^{n+1}$. Thus
\[\cycdn \circ(- \cdot X\times F^{n+1}_{1,0}) (Z) = \iota^*_{\ol{X}\times F^{n+1}} ({\rm cyc}[Z]) =  \iota^*_{\ol{X}\times \ol{F}^{n+1}_{1,0}}({\rm cyc}[\ol{Z}]) =  \iota^*_{\ol{X}\times \ol{F}^{n+1}_{1,0}} ( \beta ( \ol{Z})) \]
in \[K_0^{|W'|}(\ol{X}\times F^{n+1}_{1,0} / \ol{X}\times F^{n})^{(d+n)} \cong K_0^{|W'|}((\ol{X};D)\times \ol{\square}^n; (\ol{X};D)\times \partial \ol{\square}^n / (\ol{X};D)\times F^n)^{(d+n)}.\] 
where the isomorphism follows from the fact that $W'$ is disjoint from $D\times \cub{n+1}$.
        Let $K_0^{d+n}(\ol{X}|D, n)$ be the direct limit 
        \[K_0^{d+n}(\ol{X}|D, n) = \varcolim_{P\in (z^{d+n}(\ol{X}|D, n))_\Q} K_0^{|P|}((\ol{X};D)\times \ol{\square}^{n}; (\ol{X};D)\times \partial  \ol{\square}^n / (\ol{X};D)\times F^n)^{(d+n)}. \]
        By construction, the cycle class map ${\cycdn}$ factors through $K_0^{d+n}(\ol{X}|D, n)$, and we have a commutative diagram
        \[ \xymatrix{  N(z^{d+n}(\ol{X}|D, n+1))_\Q \ar[r]^-{\beta_{d+n}} \ar[d]^{\cdot X\times F^{n+1}_{1,0}}&  K_0^{d+n}(\ol{X}|D, n+1) \ar@/^1pc/[rd]^{\Psi} \ar[d]^{\iota^*_{\ol{X}\times F^{n+1}_{1,0}}} & \\ 
          z^{d+n}(\ol{X}|D, n)_{\Q} \ar@/_2pc/[rr]_{\cycdn} \ar[r] & K_0^{d+n}(\ol{X}|D, n) \ar[r]& K_0((\ol{X};D)\times \ol{\square}^n; (\ol{X};D)\times \partial \ol{\square}^n / (\ol{X};D)\times F^n)^{(d+n)}
        }\]        
        % \[ \xymatrix{  N(z^{d+n}(\ol{X}|D, n+1))_\Q \ar[r]^-{\beta_{d+n}} \ar[d]^{\cdot X\times F^{n+1}_{1,0}}&  K_0^{d+n}(\ol{X}|D, n+1) \ar@/^2pc/[dd]^{\Psi} \ar[d]^{\iota^*_{\ol{X}\times F^{n+1}_{1,0}}}  \\
        %   z^{d+n}(\ol{X}|D, n)_{\Q} \ar@/_1pc/[rd]_{\cycdn} \ar[r] & K_0^{d+n}(\ol{X}|D, n) \ar[d] \\
        %   & K_0((\ol{X};D)\times \ol{\square}^n; (\ol{X};D)\times \partial \ol{\square}^n / (\ol{X};D)\times F^n)^{(d+n)}
        % }\]
Now note that the map $\Psi$ has to factor through
\[{K_0}((\ol{X};D)\times \ol{\square}^{n+1}; (\ol{X};D)\times \partial'\ol{\square}^{n+1}/(\ol{X};D)\times F^{n+1})^{(d+n)}\]
But the latter group is zero. In fact, let $Y$ be either $\ol{X}\times \ol{\square}^{n}$ or $D\times \ol{\square}^n$ and identify $Y\times \P^1$ with $\ol{X}\times \ol{\square}^{n+1}$ or with $D\times \ol{\square}^{n+1}$ accordingly. We have then the exact sequence of relative $K$-groups
\[\ldots\to K_p( Y\times \P^1; Y\times \{0\}/ Y\times \{1\}) \to K_p(Y\times \P^1 /  Y\times \{1\} ) \to K_p(Y\times \{0\}) \to \ldots
\]
However, the projective bundle formula tells us that the maps $K_p(Y\times \P^1 /  Y\times \{1\} ) \to K_p(Y\times \{0\})$ are all isomorphisms, so that the groups $K_p( Y\times \P^1; Y\times \{ 0\}/ Y\times \{1\})$ are trivial. Thus the composition $\cycdn \circ(- \cdot X\times F^{n+1}_{1,0})$ is trivial, proving the Theorem.
%           The above discussion gives the following commutative diagram
%         \begin{equation}\label{eq:key-diagram-cycle-class-map}\xymatrixcolsep{3pc} \xymatrix{
%         N(z^{d+n}(\ol{X}|D, n+1))_\Q \ar[r]^-{\beta_{d+n}} \ar[d]^{\cdot X\times F^{n+1}_{1,0}}& \widetilde{K_0}((\ol{X};D)\times \ol{\square}^{n+1}; (\ol{X};D)\times \partial'\ol{\square}^{n+1}/(\ol{X};D)\times F^{n+1})^{(d+n)}\ar[d]^{\iota^*_{\ol{X}\times F^{n+1}_{1,0}}}\\
%         z^{d+n}(\ol{X}|D, n)_{\Q} \ar[r]^-{\cycdn}& K_0((\ol{X};D)\times \ol{\square}^n; (\ol{X};D)\times \partial \ol{\square}^n / (\ol{X};D)\times F^n)^{(d+n)}.
%         }\end{equation}
% Let $Y$ be either $\ol{X}\times \ol{\square}^{n}$ or $D\times \ol{\square}^n$ and identify $Y\times \P^1$ with $\ol{X}\times \ol{\square}^{n+1}$ or with $D\times \ol{\square}^{n+1}$ accordingly. We have then the exact sequence of relative $K$-groups
% \[\ldots\to K_p( Y\times \P^1; Y\times \{\infty\}/ Y\times \{1\}) \to K_p(Y\times \P^1 /  Y\times \{1\} ) \to K_p(Y\times \{0\}) \to \ldots
% \]
% However, the projective bundle formula tells us that the maps $K_p(Y\times \P^1 /  Y\times \{1\} ) \to K_p(Y\times \{0\})$ are all isomorphisms, so that the groups $K_p( Y\times \P^1; Y\times \{\infty\}/ Y\times \{1\})$ are trivial. In particular, it is trivial the group in the top right corner of \eqref{eq:key-diagram-cycle-class-map}. Thus the composition $\cycdn \circ(- \cdot X\times F^{n+1}_{1,0})$ is trivial, proving the Theorem.
        \end{proof}
\end{thm}
The following remark comes from a discussion with Kay R\"ulling.
\begin{rmk} In \cite{RulSaito}, R\"ulling and Saito constructed regulator maps 
    \[\phi^{r}_{\ol{X}|D, {\rm Nis}}\colon \mathbb{Z}(r)_{\ol{X}|D, {\rm Nis}}\to \tau_{\geq r}\mathbb{Z}(r)_{\ol{X}|D, {\rm Nis}}\to \mathcal{K}^{M}_{r, \ol{X}|D, {\rm Nis}}[-r] \]
    in the bounded derived category $D^b(\ol{X}_{\rm Nis})$, when $D_{\rm red}$ is a strict normal crossing divisor, and they showed that they induce isomorphisms \[\phi^{r}_{\ol{X}|D, {\rm Nis}}\colon H^{d+r}_{\mathcal{M}, {\rm Nis}}(\ol{X}|D, \mathbb{Z}(r)) \xrightarrow{\simeq} H^{d}(\ol{X}_{\rm Nis}, \mathcal{K}^M_{r, \ol{X}|D, {\rm Nis}}).\]
    Here, $\mathcal{K}^M_{r, \ol{X}|D, {\rm Nis}}$ denotes the Nisnevich sheaf of relative Milnor $K$-theory (its definition is rather subtle, see \cite[2.3]{RulSaito}, which is different from \cite{MR862639}). For $r=d +n$, we get by composition a natural map
    \begin{equation}\label{eq:Rulling-Saito}\CH^{d+n}(\ol{X}|D, n)_{M_{\text{ssup}}} \to  H^{2d+n}_{\mathcal{M}, {\rm Nis}}(\ol{X}|D, \mathbb{Z}(d+n)) \xrightarrow{\simeq} H^{d}(\ol{X}_{\rm Nis}, \mathcal{K}^M_{d+n, \ol{X}|D, {\rm Nis}}).\end{equation}
    If $D=\emptyset$, the natural map $\mathcal{K}^M_{d+n, \ol{X}, {\rm Nis}} \to \mathcal{K}_{d+n, \ol{X}, {\rm Nis}}$, composed with the edge homomorphism of the Nisnevich descent spectral sequence gives a further map to $K_n(\ol{X})$, and the total map from the higher Chow groups is known to agree with the cycle class map defined by Bloch and Levine. It is therefore natural to ask if there is a similar map in the relative setting 
     \begin{equation}\label{eq:dash-map}   \CH^{d+n}(\ol{X}|D, n)_{M_{\text{ssup}}} \xrightarrow{\eqref{eq:Rulling-Saito}} H^{d}(\ol{X}_{\rm Nis}, \mathcal{K}^M_{d+n, \ol{X}|D, {\rm Nis}}) \dashrightarrow  K_n(\ol{X}, D) \end{equation}
 and if this agrees with the cycle class map which we consider in this work (after tensoring with $\mathbb{Q}$). For example, if $n=0$ the map \eqref{eq:Rulling-Saito} is known to be an isomorphism when $k$ is a finite field after taking the limit over all divisors $D$ supported on $\ol{X}\setminus |D|$ as consequence of \cite{KS1} and \cite{MR862639}, and it's clear in this case that the composite map to (the limit of) $K_n(\ol{X};D)$ agrees with (the limit of) ${\rm cyc}^d_{\ol{X}|D}$ (see \cite[Theorem 1.8]{BK} for a result with fixed $D$). Note that before taking the limit, the existence of the dashed map in \eqref{eq:dash-map} might depend on the definition of relative Milnor $K$-theory. Taking the limit over all thickenings of $D$, one can apply \cite[Proposition 2.8]{RulSaito}, making the distinction between the Milnor $K$-theory defined in \cite{MR862639} and in \cite{RulSaito} irrelevant.
\end{rmk}
\begin{rmk}\label{rmk:charp}As the referee pointed out, when ${\rm char}(k) =p>0$, a Frobenius argument can be used to show that there is an isomorphism
	\begin{equation}\label{eq:indipe} \CH^{d+n}(\ol{X}|D, n)_\Q \xrightarrow{\cong} \CH^{d+n}(\ol{X}|D_{\rm red}, n)_\Q.\end{equation}
This is a particular case (in the $M_{\rm ssup}$ setting) of \cite[Theorem 1.3]{Miyazaki}, i.e.~the independence of the higher Chow groups with modulus from  the multiplicity of the divisor $D$ in positive characteristic and after inverting $p$. This is in turn a special instance of a general type of phenomena in characteristic $p$, see \cite[Theorem 3.5]{tor-div-rec} and \cite[Corollary 4.2.6]{KSY-RecAndMotives}. 

Given \eqref{eq:indipe}, the task in the statement of Theorem \ref{thm:main-theorem-cycle-class} becomes then to construct a cycle class map to the relative $K$-group $K_n(\ol{X}; D_{\rm red})^{(d+n)}$. Since $|D_{\rm red}|$ is assumed to be a strict normal crossing divisor on $\ol{X}$, one can follow the argument in the proof of \cite[Lemma 2.2]{LevineBlochrevisited} to further show the vanishing of $K_1^{|\ol{Z}|}(D\times \cub {n+1}  / D\times F^{n+1}_{j,1})^{(d+n)}$. Indeed, we can remove, without altering this group, all points of intersection, if any, between $\ol{Z}$ and $D_i\cap D_j$ for $i\neq j$ different regular components of $D$, thus reducing to the case where $D$ is regular that allows us to apply directly \cite[(2.1)]{LevineBlochrevisited}. In particular, in this situation the lifting $\beta(\ol{Z})$ constructed in \eqref{eq:fundamental-class-beta-with-support} is unique. This gives a simplification of the proof of Theorem \ref{thm:main-theorem-cycle-class} in positive characteristic. Of course, in characteristic zero the isomorphism \eqref{eq:indipe} is far from being true (the reader can convince herself by contemplating the case $n=0$), and to prove our main result we need anyway the detour of Section \ref{sec:notation-classes-to-lift-2}. 
	\end{rmk}
\subsection{Relationship with $K$-groups of geometric type}\label{ssec:Cor-Sugiyama} Before completing the proof of Claim \ref{main-claim-vanishing-class}, we explain how the cycle class map can be used to relate higher relative $K$-groups with certain $K$-groups of geometric type. Assume in this subsection that the base field $k$ is perfect.

In their work towards the construction of a theory of motives without $\mathbb{A}^1$-invariance, Kahn, Saito and Yamazaki introduced in \cite{KSY} the notion of \textit{reciprocity (pre)sheaf}, as generalization of Voevodsky's notion of homotopy invariant (pre)sheaf with transfers. In \cite{KSY2} and \cite{KSY-RecAndMotives}, the authors proposed an alternative, stronger, definition, which goes under the name of \textit{SC-reciprocity property} (where SC stands for ``Suslin-complex'') for (pre)sheaves, and that fits well with the categorical framework of the theory of motives with modulus. 

In particular, for $M = (\ol{X};D)$ a proper modulus pair (i.e., a pair consisting of a proper $k$-scheme $\ol{X}$ and an effective Cartier divisor $D$ such that $X = \ol{X}\setminus D$ is smooth over $k$), there is a distinguished SC-reciprocity presheaf with transfers $h_0(M)$ such that $h_0(M)(\Spec{k}) = \CH_0(\ol{X}|D)$ (see \cite[Section 7]{KSY2}). The subcategory $\mathbf{RSC}\subset \mathbf{PST}$ of SC-reciprocity presheaves   comes equipped with a tensor structure, unlike the category of reciprocity (pre)sheaves.

\subsubsection{}In \cite{Sugiyama:2017aa}, Sugiyama introduced a notion of $K$-group of geometric type $K^{\rm geo}(-)$, analogue in the modulus setting to the $K$-group of geometric type of Kahn-Yamazaki \cite{KahnYamazakiDuke}, with the goal of giving a description in terms of symbols of the tensor product of SC-reciprocity presheaves. In particular, he obtained
\begin{thm}[Theorem 3.1 \cite{Sugiyama:2017aa}]Let $F_1,\ldots, F_n$ be SC-reciprocity presheaves with transfers. Then there is an isomorphism 
    \[K_{\rm sum}^{\rm geo}(k;F_1,\ldots, F_n)\cong (F_1\tensor_{\RSC} F_2\tensor_{\RSC}\tensor\ldots \tensor_{\RSC} F_n)(\Spec{k}).\]
    \end{thm}
    The ``sum'' refers, as in our Definition of Section \ref{sec:moduli-conditions}, to a chosen modulus condition among symbols (an alternative notion, involving a ``max''-modulus condition, is also discussed). In some special cases, the $K$-groups $K^{\rm geo}(-)$ are related to algebraic cycles, and this is what we exploit in the following Corollary to Theorem \ref{thm:main-theorem-cycle-class}.
\begin{cor}\label{cor:Sugiyama-Cor} Let $\ol{X}$ be a smooth, proper $k$-variety of dimension $d$ and let $D$ be an effective Cartier divisor on $\ol{X}$ such that $|D_{\rm red}|$ is strict normal crossing. Then there is a canonical homomorphism
   \[{\rm cyc}_{\ol{X}|D, n}^{{\rm geo}} \colon K_{\rm sum}^{\rm geo}(k;h_0(\ol{X};D), \mathbb{G}_m, \mathbb{G}_m, \ldots, \mathbb{G}_m)\tensor\mathbb{Q}\to K_n(\ol{X};D)^{(d+n)}\]
    for each $n\geq 0$.
    \begin{proof}By \cite[Proposition 5.7]{Sugiyama:2017aa}, there is a surjection 
        \begin{equation}\label{eq:cycRin}K_{\rm sum}^{\rm geo}(k;h_0(\ol{X};D), \mathbb{G}_m, \mathbb{G}_m, \ldots, \mathbb{G}_m) \to \CH^{d+n}(\ol{X}|D, n)_{M_\Sigma}\end{equation}
        (where the $K$-group on the left hand side has $n$-copies of $\mathbb{G}_m$). It is however clear from the proof that the map \eqref{eq:cycRin} factors through $\CH^{d+n}(\ol{X}|D, n)_{M_{\rm ssup}}$. We recall Sugiyama's argument for the reader's convenience. The group $K_{\rm sum}^{\rm geo}(k; h_0(\ol{X};D),  \mathbb{G}_m, \ldots, \mathbb{G}_m)$ is a quotient of
        \[(h_0(\ol{X};D)\tensor^M \mathbb{G}_m\tensor\ldots \tensor^M \mathbb{G}_m)(k),\]
         where $\tensor^M$ denotes the tensor product of Mackey functors (see \cite[2.8]{KahnYamazakiDuke}).  A typical relation in $K_{\rm sum}^{\rm geo}(-)$ as above is of the form $\sum_{c\in C^o}v_c(f){\rm Tr}_{k(c)/k}(\alpha(c)\tensor g_1(c)\tensor \ldots g_n(c))$, where $C$ is a proper integral curve over $k$, $C^o$ is the open complement of a modulus $\mathfrak{m}$ on $C$, $f$ is a function on $C$ having modulus $\mathfrak{m}$, $g_i$ are invertible functions on $C^o$ and $\alpha$ is an elementary admissible finite correspondence from $(C, \mathfrak{m})$  to $(\ol{X};D)$ (see \cite[Definition 3.2]{Sugiyama:2017aa} for details). If $\ol{Z}$ denotes the closure of $\alpha$ in $C\times \ol{X}$, the above relation gives a $1$-cycle in $\ol{X}\times \ol{\square}^{n+1}$ image of 
        \[\phi :=q\times p^*f\times p^*g_1\ldots \times p^* g_n \colon \ol{Z}^N \to \ol{X}\times  \ol{\square}^{n+1} \]
        where $p$ and $q$ denote the two projections to $C$ and $\ol{X}$. Then one has
        \[\phi^*(   D\times \ol{\square}^{n+1}) \leq \phi^* (\ol{X}\times \{1\}\times \ol{\square}^{n})\]
        so that the cycle $W :={\rm Im}(\phi)\cap  X\times \square^{n+1}$ satisfies the $M_{\rm ssup}$-modulus condition with respect to the face $\ol{X}\times F^{n+1}_{1,1}$. In particular, it is a relation in $\CH^{d+n}(\ol{X}|D,n)_{M_{\rm ssup}}$ and we get a (surjective) map
        \begin{equation}\label{eq:cycRin-Msup}K_{\rm sum}^{\rm geo}(k; h_0(\ol{X};D), \mathbb{G}_m, \mathbb{G}_m, \ldots, \mathbb{G}_m) \to \CH^{d+n}(\ol{X}|D, n)_{M_{\rm ssup}}.\end{equation}
       Composing \eqref{eq:cycRin-Msup} with ${\rm cyc}^{d+n}_{\ol{X}|D}$ gives the required homomorphism ${\rm cyc}^{\rm geo}_{\ol{X}|D, n}$.
        \end{proof}
    \end{cor}
\begin{rmk}The same argument can be used to construct a cycle class map from the variant $K^{\rm geo}_{\rm max}(-)$. We leave the details to the reader.
	\end{rmk}
\subsection{Lifting classes}\label{sec:lifting-classes} In this Section we finally explain how to construct the desired lifting. We resume the notations of \ref{sec:notation-classes-to-lift}-\ref{sec:notation-classes-to-lift-2}, so suppose that $Z$ is a cycle in $N(z^{d+n}(\ol{X}|D, n+1))$, with irreducible components $Z_1,\ldots Z_r$. %Note that the class ${\rm cyc}_{\ol{Z}}[\ol{Z}_k]\in K_0^{|\ol{Z}|}(\ol{X}\times \cub {n+1})^{(d+n)}$ is the image of ${\rm cyc}_{\ol{Z}_k}[\ol{Z}_k] \in K_0^{|\ol{Z}_k|}(\ol{X}\times \cub {n+1})^{(d+n)}$ via the natural map
%\[\rho_{\ol{Z}_k, \ol{Z}}\colon  K_0^{|\ol{Z}_k|}(\ol{X}\times \cub {n+1})^{(d+n)} \to K_0^{|\ol{Z}|}(\ol{X}\times \cub {n+1})^{(d+n)}\]
%In particular, the class ${\rm cyc}_{\ol{Z}}[\ol{Z}]$ is the sum $\sum_{k=1}^r \rho_{\ol{Z}_k, \ol{Z}} {\rm cyc}_{\ol{Z}_k}[\ol{Z}_k]$.
 As explained in the previous Section, it is enough to show that each ${\rm cyc}_{\ol{Z}_k}[\ol{Z}_k]$ projects to a class in $K_0^{|\Zbk|}(\ol{X}\times \cub {n+1}  / \ol{X}\times F^{n+1}_{j,1})^{(d+n)}$ that vanishes when restricted to $K_0^{|\Zbk|}(D\times \cub {n+1}  / D\times F^{n+1}_{j,1})^{(d+n)}$ (for some $j$ depending on $k$). In order to show it, we exploit the modulus condition.

By definition, the modulus condition on a cycle $W$ is tested on its irreducible components. Consider now the case of our cycle $Z$. The strong sup-modulus conditions, satisfied by each $Z_k$, takes the following form. Write $\ol{Z}_k$ for the closure of $Z_k$ in $\ol{X}\times \ol{\square}^{n+1}$.  Let $\phi_{\ol{Z}_k}\colon \ol{Z}^N_k\to \ol{X}\times \ol{\square}^{n+1}$ be the normalization morphism followed by the natural inclusion. Then, there exists $j = j(k)\in \{1,\ldots, n+1\}$ such that 
\[ \phi_{\ol{Z}_k}^* (D\times \ol{\square}^{n+1} ) \leq \phi_{\ol{Z}_k}^* (\ol{X}\times F^{n+1}_{j,1}).\]

In particular, we have an inclusion of sets $\ol{Z}_k\cap D\times \ol{\square}^{n+1} \subseteq \ol{X}\times F^{n+1}_{j,1}$. \
%Without loss of generality, we can assume that the intersection of $\ol{Z}_k$ with the divisor $D\times \ol{\square}^{n+1}$ is given by a single closed point $P$, that coincides with the intersection of $\ol{Z}_k$ with $\ol{X}\times  F^{n+1}_{j,1}$. 
For any point $P\in \ol{Z_k} \cap D\times \ol{\square}^{n+1}$, we have the following commutative diagram
  %K_0^{|\Zbk|}((\ol{X}, D)\times \cub {n+1} / \ol{X}\times F^{n+1}_{j,1})^{(d+n)} \ar[r] 
%K_0^{|\Zbk|}((\ol{X}, D)\times \cub {n+1})^{(d+n)} \ar[u]
\begin{equation}\label{eq:fund-diagram-mod-cond-gives-lifting} \xymatrix{K_0^{|\Zbk|}(\ol{X}\times \cub {n+1}  / \ol{X}\times F^{n+1}_{j,1})^{(d+n)} \ar[r] & K_0^{|P|}(D\times \cub {n+1}  / D\times F^{n+1}_{j,1})^{(d+n)}\\
 K_0^{|\Zbk|}(\ol{X}\times \cub {n+1})^{(d+n)}\ar[u]^{p^j_{\ol{X}}} \ar[r]  & K_0^{|P|}(D\times \cub {n+1})^{(d+n)}\ar[u]\\
 K_0^{|P|}(\ol{X}\times F^{n+1}_{j,1})^{(d+n-1)}\ar[u]^{i^{n+1}_{j,1, *}} \ar[r] & K_0^{|P|}(D\times F^{n+1}_{j,1})^{(d+n-1)} \ar[u]^{i_{D\times F^{n+1}_{j,1}, *}}
}
\end{equation}
in which the middle horizontal map is induced by $\iota^*_{D\times \ol{\square}^{n+1}}$ composed with the projection 
\[ K_0^{|\Zbk \cap (D\times \cub{n+1})|}( D \times \cub {n+1})^{(d+n)} = \bigoplus_{Q\in \Zbk \cap (D\times \cub{n+1})} K_0^{|Q|}(D\times \cub {n+1})^{(d+n)} \to K_0^{|P|}(D\times \cub {n+1})^{(d+n)}. \] 
and the top horizontal arrow is induced in a similar way by $\iota^*_{D\times \ol{\square}^{n+1}}$ composed with the projection to the $P$-component on the cofiber group. We have to show that the restriction to the $K_0$-groups with support in $P$ of $\iota_{D\times \cub {n+1}}^*([\cO_{\ol{Z}_k}]) = \iota_{D\times \cub {n+1}}^* ({\rm cyc}_{\ol{Z}_k}[\ol{Z}_k])$   is the image of a class $\alpha_{Z_k,j}$ in  $K_0^{|P|}(D\times F^{n+1}_{j,1})^{(d+n)}$ along the push-forward $i_{D\times F^{n+1}_{j,1}, *}$ for every $P$. Since we check this for every point $P$, we can assume that the intersection of $\ol{Z}_k$ with the divisor $D\times \ol{\square}^{n+1}$ is given by a single closed point, that coincides with the intersection of $\ol{Z}_k$ with $\ol{X}\times  F^{n+1}_{j,1}$.
\subsubsection{}In order to treat uniformly  the case where $\ol{Z}_k$ is not regular in a neighborhood of $D$, we change the notation a bit and consider the following slightly more general situation (see also \cite[Section 5]{BS}, where we use the same convention). Let $Y$ be a smooth (connected) $k$-variety of dimension $d+1$, equipped with a smooth divisor $F$ and an effective Cartier divisor $D$. Asume that $F$ and $D$ satisfy together the following condition:
\begin{enumerate}
\item[$(\bigstar)$]
There is no common component of $D$ and $F$, and 
$D_{red}+F$ is a (reduced) simple normal crossing divisor on $Y$.
\end{enumerate}
\begin{df}Let $C$ be an integral curve contained in $X = Y - (F+D)$. Write $\ol{C}$ for the closure of $C$ in $Y$ and $\ol{C}^N$ for the normalization of $\ol{C}$. Let $\phi_{\ol{C}}\colon \ol{C}^N\to Y$ be the natural map. We say that $C$ satisfies the \textit{modulus condition with respect to the divisor $D$ and the face $F$} if the following inequality of Cartier divisors on $\ol{C}^N$ holds:
    \[ \phi_{\ol{C}}^*(D)\leq \phi_{\ol{C}}^{*}(F)\]
 \end{df}
Write $\iota_D\colon D\to Y$ (resp. $\iota_F\colon F\to Y$) for the inclusion of $D$ (resp. of $F$) in $Y$ and write $j_{D,F}\colon D\cap F \to D$ for the inclusion of the intersection of $D$ and $F$ inside $D$.
Let $C$ be an integral curve satisfying the modulus condition. Assume (using the same argument discussed above) that $\ol{C}\cap D$ is given by a single point $P$. In the current setting, the diagram \eqref{eq:fund-diagram-mod-cond-gives-lifting} takes the following form:
\[\xymatrix{ K_0^{|\ol{C}|}(Y/F)^{(d)} \ar[r] & K_0^{|P|}(D / D\cap F)^{(d)} \\
K_0^{\ol{C}}(Y)^{(d)} \ar[u]\ar[r]^{\iota_D^*} & K_0^{|P|}(D)^{(d)}\ar[u] \\
K_0^{|P|}(F)^{(d-1)}\ar[u]^{\iota_{F,*}} \ar[r] & K_0^{|P|}(D\cap F)^{(d-1)}\ar[u]^{j_{D,F,*}}
}
\]
We write $[\ol{C}]$ for the fundamental class of $\ol{C}$ in $K_0^{\ol{C}}(Y)^{(d)}$ given by the class of the structure sheaf $\cO_{\ol{C}}$. 
\begin{prop}\label{prop:regular-case-curves-gensetting}The restriction  $\iota_D^* ([\ol{C}]) = \iota_D^*([\cO_{\ol{C}}])$ of the fundamental class of $\ol{C}$ along the divisor $D$ is the image of a class $\alpha_{\ol{C}}$ in  $K_0^{|P|}(D\cap F)$ along the push-forward ${j_{D,F,*}}$.
    \begin{proof}We start by assuming that $\ol{C}$ is regular in a neighborhood of $P$. Since $\ol{C}$ is not contained in $D$, the module $\cO_{\ol{C}}$ is $\cO(-D)$-torsion free and we have an equality 
        \[ \iota_D^* ([\cO_{\ol{C}}]) = [\cO_{\ol{C}} \tensor_{\cO_Y} \cO_D]\]
        in $K_0^{|P|}(D)$, and $\cO_{\ol{C}} \tensor_{\cO_Y} \cO_D$ is a module of finite homological dimension over $\cO_D$. The class $\iota_D^* ([\cO_{\ol{C}}])$ is supported on $P$ by assumption, and we can work locally around $P$ in the following sense. Let $\cO_P$ be the local ring of $\cO_Y$ at $P$. Since $Y$ is regular at $P$, $\cO_P$ is a Noetherian regular equicharacteristic local ring (containing the field $k$). By Cohen's structure theorem (together with the fact that $D_{\rm red}+F$ is strict normal crossing), its completion $\widehat{\cO}_P$ is then isomorphic to a power series ring $K \llbracket x_1,\ldots, x_d, t_{1} \rrbracket$, where $t_1$ is the image in $\cO_P$ of a local parameter for the smooth divisor $F$. Up to reordering, the ideal $I_D$ of the divisor $D$ in $\widehat{\cO}_P$ will be then generated by an element $\prod_{i=1}^s x_i^{m_i}$. 
        By \cite[Proposition 3.19]{TT}, we can replace the group $K_0^{|P|}(D)$ with $K_0^{|P|}(\Spec{\widehat{\cO}_P/I_D})$ and the group $K_0^{|P|}(D\cap F)$ with $K_0^{|P|}(\Spec{\widehat{\cO}_P/(I_D, t_1)})$. Let $\widehat{\cO}_{\ol{C},P}$ be the completion of the local ring of $\ol{C}$ at $P$. After this reduction, the class $\iota_D^* ([\cO_{\ol{C}}])$ we are after is  (the class of) the module  $\widehat{\cO}_{\ol{C},P} \tensor \widehat{\cO}_P/I_D$. Since $\ol{C}$ is regular at $P$, we can assume that the image of one of the paramenters $x_i$ or $t_1$  is a local parameter for $\ol{C}$. Without loss of generality (the proof is substantially identical in other cases), we assume that this role is played by $x_d$. Thus, we can write in $\widehat{\cO}_{\ol{C}, P}$ 
      \begin{equation}\label{eq:def-equ-curve-modulus-gen}  x_i = x_d^{a_i} v_i, \quad t_1 = x_d^b u_1,\quad \text{ for } i=1,\ldots, d-1.\end{equation}
         and elements $v_i, u_1 \in \widehat{\cO}_{\ol{C}, P}^{\times}$, that we can write as power series in $x_d$. Note that since the curve is actually passing through the point $P$, the exponents $a_i$ and $b$ have to be positive. The modulus condition gives then the following inequality
        \begin{equation}\label{eq:mod-cond-curve}\sum_{i=1}^{s} m_i a_i \leq b.\end{equation}
 Write $J$ for the ideal of $K \llbracket x_1,\ldots, x_d, t_{1}\rrbracket$ defined by the equations \eqref{eq:def-equ-curve-modulus-gen}: by construction, it agrees with the kernel of the map 
 \[ K \llbracket x_1,\ldots, x_d, t_{1}\rrbracket \cong \widehat{\cO}_P \to \widehat{\cO}_{\ol{C}, P}.\]
  By \eqref{eq:mod-cond-curve}, the ideal $(J, \prod_{i=1}^s x_i^{m_i})$ and the ideal $(\prod_{i=1}^s x_i^{m_i}, t_1, x_i - x_{d}^{a_i}v_i)_{i}$ coincide. Now, the module 
 \begin{align*} M_{\ol{C}, P} &= K \llbracket x_1,\ldots, x_d, t_{1} \rrbracket/(\prod_{i=1}^s x_i^{m_i}, t_1, x_i - x_{d}^{a_i}v_i)_{i}\\  
 & \cong  (K \llbracket x_1,\ldots, x_d\rrbracket/ (\prod_{i=1}^s x_i^{m_i})) \llbracket t_{1} \rrbracket/ (t_1, x_i - x_{d}^{a_i}v_i)_{i} \end{align*}
 has finite homological dimension as module over $(K \llbracket x_1,\ldots, x_d\rrbracket/ (\prod_{i=1}^s x_i^{m_i})) \llbracket t_{1} \rrbracket/(t_1)$, and is supported on $P$.
 
 It gives then a well defined class $[M_{\ol{C}, P}]$ in the $K_0$ group with support $K_0^{|P|}(\Spec{ \widehat{\cO}_{P} / (I_D, t_1)})$ which satisfies 
 \[j_{D,F,*} [M_{\ol{C}, P}] = [\widehat{\cO}_{\ol{C},P} \tensor \widehat{\cO}_P/I_D]   = \iota_D^* ([\cO_{\ol{C}}]) \]
 as required.
 
 We now deal with the case where $\ol{C}$ is not necessarily regular in a neighborhood of $\ol{C}\cap D$. Let $\phi\colon \ol{C}^N \to\ol{C}$ be the normalization morphism. It fits in a commutative diagram
\begin{equation} \label{eq:diag-final}\xymatrix{  \ol{C}^N \ar@{^{(}->}[r]^-{j}\ar[d]^{\phi} &Y \times \mathbb{P}^M \ar[d]^p \\
\ol{C} \ar@{^{(}->}[r]  & Y
}
\end{equation}
where $p$ is the natural projection. The curve $\ol{C}^N$ is now regular and embedded in the smooth variety $Y \times \mathbb{P}^M = \P^M_Y$, and satisfies the modulus condition with respect to $\P^M_D$ and the face $\P^M_F$. In particular, we can apply the normal case to conclude that the class of $\ol{C}^N$ in the cofiber group $K_0^{|\ol{C}^N|}(Y\times \P^M / F\times \P^M)$ dies when restricted to $D\times \P^M$. 
The covariant functoriality of $K$-theory for proper maps of finite Tor-dimension, as in \cite[3.16.4]{TT}, applied to the diagram \eqref{eq:diag-final}, gives an equality of $K_0$-classes $p_*[ \cO_{\ol{C}^N}] = [p_* \cO_{\ol{C}^N}] =  [ \cO_{\ol{C}}] +  [S]$, where $S$ is a coherent sheaf on $\ol{C}$ supported on finitely many points $y_1,\ldots, y_r$. The class $[p_* \cO_{\ol{C}^N}]$ in $K_0^{|\ol{C}|}(Y)^{(d)}$ maps then to a class in the cofiber group $K_0^{|\ol{C}|}(Y/ F)^{(d)}$ that is mapped to zero by construction when restricted to $K_0^{|P|}(D/ F\cap D)^{(d)}$. Let now $T\subset \ol{C}$ be  the subset of $\ol{C}$ given by the union of the (closed) points $y_i$. Since $Y$ is regular, we have an isomorphism $K_0^{|\ol{C}|}(Y)^{(d)} \xrightarrow{\cong} K_0^{|\ol{C}\setminus T|}(Y\setminus T)^{(d)}$. In particular, as the sheaf $S$ is supported on $T$, we have $p_*[ \cO_{\ol{C}^N}] = [p_* \cO_{\ol{C}^N}] =  [ \cO_{\ol{C}}]$ in $K_0^{|\ol{C}|}(Y)^{(d)}$. We can thus replace $ [ \cO_{\ol{C}}]$ with the push-forward class $[p_* \cO_{\ol{C}^N}]$ and we are done.
        \end{proof}
    \end{prop}
Applying Proposition \ref{prop:regular-case-curves-gensetting} to the setting of \ref{sec:lifting-classes}, we can deduce the following Proposition, proving Claim \ref{main-claim-vanishing-class}. 
\begin{prop}\label{prop:regular-case-curves}The restriction $\iota_{D\times \cub {n+1}}^*([\cO_{\ol{Z}_k}]) = \iota_{D\times \cub {n+1}}^* ({\rm cyc}_{\ol{Z}_k}[\ol{Z}_k])$ of the fundamental class of $\ol{Z}_k$ along the divisor $D\times \ol{\square}^{n+1}$ is the image of a class $\alpha_{Z_k,j}$ in  $K_0^{|P|}(D\times F^{n+1}_{j,1})^{(d+n)}$ along the push-forward $i_{D\times F^{n+1}_{j,1}, *}$.
       \end{prop}

\noindent\emph{Acknowledgments.} The core of this paper was written in 2014, while I was a graduate student in Essen under the supervision of Marc Levine. It is a pleasure to thank him heartily for much advice, constant support and encouragement. I would also like to thank Kay R\"ulling for valuable comments on a preliminary version of this paper. Finally, I would like to thank the referee for a very careful reading of the manuscript and for suggesting many improvements.

    \bibliography{bibMilnorChowMod} 
    \bibliographystyle{siam}
%\printindex
\end{document}